\documentclass{amsart}
\usepackage{hyperref,bbm,booktabs,graphicx,subcaption,enumitem,wrapfig}
\usepackage[final,expansion=true,protrusion=true,spacing=true,kerning=true]{microtype}
\usepackage[font=footnotesize]{caption}
\usepackage[ 
    hmarginratio={1:1},     
    vmarginratio={1:1},     
    heightrounded,          
]{geometry}
\title{Fast high-dimensional node generation with variable density}
\author{O.  Vlasiuk}
\author{T. Michaels}
\author{N. Flyer}
\author{B. Fornberg}
\date{\today}
\address{Center for Constructive Approximation, Department of Mathematics, Vanderbilt University, Nashville, TN 37240, USA}
\email{oleksandr.vlasiuk@vanderbilt.edu}
\email{timothy.j.michaels@vanderbilt.edu}
\address{Analytics and Integrated Machine Learning, National Center for Atmospheric Research, Boulder, CO 80305, USA}
\email{flyer@ucar.edu}
\address{Department of Applied Mathematics, University of Colorado, Boulder, CO 80309, USA}
\email{bengt.fornberg@colorado.edu}

\newcommand{\fwdth}{0.45} 
\newcommand{\bgwdth}{0.8} 
\newcommand{\smwdth}{0.37}

\newcommand{\bs}[1]{\boldsymbol{#1}} 

\newcommand{\ndv}{\mathcal{V}} 
\newcommand{\lat}{\mathcal{L}}
\newcommand{\rmin}{\mathcal{M}}
\newcommand{\Fr}{\mathrm{frac}}
\newcommand{\VoxelsPerSide}{M}
\newcommand{\tabbox}{\parbox[c][][c]{.5\textwidth}}
\def\tabline{\\ } 
\newcommand{\OffsetLattice}{C_1}
\newcommand{\OffsetGrad}{{C}_2}
\newcommand{\ScalingFactor}{f}
\newcommand{\ShiftVector}{\boldsymbol{h}}
\newcommand{\VarRadius}{r}
\newcommand{\VarAngle}{a}
\newcommand{\VarPolarAngle}{p}
\newcommand{\AngleSpacing}{B}

\newcommand{\nd}{\boldsymbol{x}}
\newcommand{\vecx}{\boldsymbol{x}}

\newlist{algo}{enumerate}{1} 
\setlist[algo, 1]{label = \texttt{Step \arabic*}, start=0, align=left,  labelwidth=0pt, itemindent=!} 
\newlist{algo2}{enumerate}{1} 
\setlist[algo2, 1]{label = \texttt{Step \arabic*}, start=3, align=left,  labelwidth=0pt, itemindent=!} 
\newlist{algo3}{enumerate}{1} 
\setlist[algo3, 1]{label = \texttt{Step \arabic*}, start=5, align=left,  labelwidth=0pt, itemindent=!} 

\makeatletter
\def\subsubsection{\@startsection{subsubsection}{3}%
  \z@{.5\linespacing\@plus.7\linespacing}{-.5em}%
  {\normalfont\bfseries}}
\makeatother

\begin{document} 
\maketitle 
\begin{abstract}
    We present an algorithm for producing discrete distributions with a prescribed nearest-neighbor distance function. Our approach is a combination of quasi-Monte Carlo (Q-MC) methods and weighted Riesz energy minimization: the initial distribution is a stratified Q-MC sequence with some modifications; a suitable energy functional on the configuration space is then minimized to ensure local regularity. The resulting node sets are good candidates for building meshless solvers and interpolants, as well as for other purposes where a point cloud with a controlled separation-covering ratio is required. Applications of a three-dimensional implementation of the algorithm, in particular to atmospheric modeling, are also given.
\end{abstract}
\section{Introduction} 
\subsection{RBFs and meshless methods}
In a number of important applications, usefulness of meshless methods in general, and of
radial basis functions (RBFs) in particular, is well-known. They have found
their way into high-dimensional interpolation, machine learning, spectral
methods, vector-valued approximation and interpolation, just to name a few
\cite{Wendland2004, Chang2011, Fornberg1996, BroomheadD.S.andLowe1988, Shankar2014}.  RBFs have multiple advantages, most importantly extreme flexibility in forming
stencils (in the case of RBF-FD) and high local adaptivity; allowing spectral
accuracy on irregular domains; the fact that the corresponding interpolation
matrix (denoted by $ \bs{A} $ below) is positive definite for several types of
radial functions and does not suffer from instability phenomena characteristic
of some of the alternative interpolation methods.

Applying RBF-FD stencils to building solvers requires an efficient way of
distributing the \textit{nodes} of basis elements in the domain, which can be either a
solid or a surface. The tasks of modeling and simulation often call for
massive numbers of nodes, so it is important to ensure that the distribution process is
easily scalable. One further has to be able to place the RBFs according to a
certain density, as a method of local refinement, for example, at the boundary,
or in regions of special interest. Yet another challenge arises when it is
necessary to deal with complex or non-smooth domains and/or surfaces.

Recall \cite{Fornberg2015b} that an RBF is a linear combination of the form
\begin{equation}\label{eq:rbf}
    \mathcal S(\vecx) = \sum_{k=1}^{K} c_k \phi(\|\vecx-\vecx_k\|), 
\end{equation} where $ \phi(\cdot) $ is a radial
function, and $ \vecx_k,\ k=1,\,\ldots,K $, is a collection of pairwise distinct
points in $ \mathbb{R}^d $. A common choice of $ \phi $ is the Gaussian $ \phi(r) = e^{-(\epsilon r)^2} $, although one may also use $ 1/(1+(\epsilon r)^2)$,  $r^{2 p}\log(r),\ p\in\mathbb N $, etc. In this discussion, we are not concerned with the
distinctions between the different radial kernels, so the reader can assume
that $ \phi(r) = e^{-(\epsilon r)^2} $. In contrast to pseudospectral methods \cite{Fornberg1996}, RBF-FD approach means that to obtain a useful approximation of a function, or a differential operator, the nodes in expressions like \eqref{eq:rbf} must be in the vicinity of the point $ \vecx $, and therefore a large number of stencils are constructed throughout the underlying set.  It is well-known that the matrix 
\begin{equation}
    \label{eq:Amat}
    \bs{A} =
    \begin{bmatrix} \phi(\|\vecx_1-\vecx_1\|)& \phi(\|\vecx_1-\vecx_2\|) &
        \ldots & \phi(\|\vecx_1-\vecx_K\|) \\ \phi(\|\vecx_2-\vecx_1\|)&
        \phi(\|\vecx_2-\vecx_2\|) & \ldots & \phi(\|\vecx_2-\vecx_K\|) \\
        \vdots& \vdots  &  & \vdots  \\ \phi(\|\vecx_K-\vecx_1\|)&
        \phi(\|\vecx_K-\vecx_2\|) & \ldots & \phi(\|\vecx_K-\vecx_K\|) \\
    \end{bmatrix} 
\end{equation}
is positive definite if the nodes $ \vecx_1\ldots \vecx_K $
are all distinct \cite{Schoenberg1938}, and so under this assumption there exist  $ K
$-point RBF interpolants for any function data. A different question, however,
is whether the matrix $ \bs{A} $ will be well-conditioned: it is not the case, for
example, when the nodes are placed on a lattice and $ \epsilon\to 0 $,
\cite{Fornberg2015}. The other extreme, having low regularity, also does not provide a reliable source of nodes, as can be seen on the example of the Halton sequence   \cite{Fornberg2015}. Furthermore, node clumping  can lead to instability of PDE solvers, \cite{Fornberg1996}. To avoid this, one must guarantee that the nodes are separated. In effect, generally the quasi-uniform node sets generated by the present algorithm, or, for example, the one constructed by the third and fourth authors \cite{Fornberg2015a},  perform better than either lattice or the Halton sequence.

In many applications, one has to ensure that the distance from a node $ \nd $
to its nearest neighbor behaves approximately as a function of the position of
the node \cite{Fornberg2015a}. Prescribing this function, $ \rho(\nd) $, which
we call the \label{p:radial_density}\textit{radial density}, is a natural way to treat the cases when a
local refinement is required in order to capture special features of the
domain. In the present paper we will describe a method of node placement for
which the actual distance to the nearest neighbor, denoted by $ \Delta(\nd) =
\min_{\nd'\neq \nd} \|\nd'-\nd\| $, satisfies the above description. To
summarize, we are interested in a procedure for obtaining discrete
configurations inside a compact set that will: 
\begin{itemize}[leftmargin=*] 
    \item
        guarantee that $ \Delta(\nd) \asymp \rho(\nd) $ (that is, differ only up to
        a constant factor) for a given function $ \rho(\nd) $ with a reasonably
        wide choice of $ \rho $; 
    \item be suitable for mesh-free PDE
        discretizations using RBFs, i.e., produce well-separated configurations
        without significant node alignment; 
    \item result in quasi-uniform node distributions also on the surface boundaries of the domain; \item be
        computationally efficient, easily scalable, and suitable for
        parallelization.  
\end{itemize}

\subsection{Notation and layout}
The bold typeface is reserved for vectors in $ \mathbb R^d $. With few exceptions, letters of the Greek alphabet denote functions, calligraphic letters and $ \Omega $ denote sets, and the regular Roman typeface is used for scalar variables. The symbolic notation employed throughout the paper is summarized in Table~\ref{tab:notation}.

The paper is structured as follows: Section \ref{subsec:rbfs} outlines the RBF-FD method using Gaussian and Polyharmonic Spline kernels;  Sections \ref{subsec:riesz} and \ref{subsec:qmc} introduce the two essential components of our approach, Riesz energy functionals and quasi-Monte Carlo methods. The main algorithm and its discussion are the subjects of Sections \ref{subsec:algo} and \ref{subsec:algo_discussed}, respectively. Sections~\ref{subsec:atm_nodes}--\ref{subsec:sph_shell} offer applications of the algorithm; the corresponding run times are summarized in Section~\ref{subsec:timings}. Section~\ref{subsec:comparisons} contains comparisons of the condition numbers of RBF-FD matrices with stencils on periodic Riesz minimizers, Halton nodes, and the Cartesian grid; Section~\ref{subsec:applicability} discusses the range of dimensions where the present method is applicable. The \hyperref[appendix]{Appendix} is dedicated to numerical experiments with the mean and minimal separation distance of Riesz minimizers and irrational lattices.

\begin{table}[t]
    \centering
    \[
        \begin{array}{cl}\toprule
            \makebox[.25\textwidth]{Symbol}  & \makebox[.45\textwidth]{Description} \\ \midrule
            \bar x                     & \tabbox{mean value of $ x $}\tabline 
            \alpha_1, \alpha_2, \ldots, \alpha_{d-1}    & \tabbox{fixed linearly independent irrational numbers}  \tabline 
            \mathcal{C}                                 & \tabbox{unit cube $ [0,1]^d $} \tabline
            \OffsetLattice, \OffsetGrad                 & \tabbox{positive constants in \eqref{eq:irr_lat} and \eqref{eq:node_update}}\tabline
            \bs{c}_m                                      & \tabbox{the closest to the origin corner of $ \mathcal V_m $} \tabline
            \chi(\cdot\,; \Omega)                         & \text{characteristic function of the set $ \Omega $}\tabline
            \mathcal{D}                                & \tabbox{$ m\in\mathcal D $ if $ \mathcal V_m $ has nonempty neighbors on \ref{step:subcubes}} \tabline
            \Delta(\vecx)                               & \tabbox{distance from the node $ \vecx $ to the nearest neighbor} \tabline 
            \Delta(\{\vecx_1,\ldots\vecx_N\})           & \tabbox{separation of the configuration } \tabline 
            \Delta^{k}(\vecx),\ k=1,\ldots,K            & \tabbox{distance from $ \vecx $ to the $k$-th nearest neighbor}\tabline 
            \Delta_n                     & \tabbox{minimal separation of periodized $ \lat_n $, $ \rmin_n $, p. \pageref{p:delta_n}}\tabline
            \bar\Delta_n                            & \tabbox{mean separation of the periodized $ \lat_n $, $ \rmin_n $, p. \pageref{p:delta_n}}\tabline
            \mathcal E                                  & \tabbox{$ m\in\mathcal{E} $ if $ \mathcal V_m $ is empty after \ref{step:densefill}} \tabline
            E_s,\ E_s^\kappa                            & \tabbox{Riesz $s$-energy and weighted $s$-energy, \eqref{eq:s-riesz}-\eqref{eq:sw-riesz}} \tabline  
            \Fr( x )                                     & \text{fractional part of the nonnegative number }x \tabline 
            \mathcal{H}_d                               & d\text{-dimensional Hausdorff measure}\tabline 
            \kappa                                      & \tabbox{kernel of the weighted Riesz energy, \eqref{eq:sw-riesz}} \tabline 
            \lambda                                     & \tabbox{interpolated inverse of  $ \bar\Delta_n $, p. \pageref{p:delta_n}}\tabline 
            \lat_n,\ n\geq1                             & n\tabbox{-point irrational lattice, \eqref{eq:irr_lat}} \tabline 
            \rmin_n,\ n\geq1                            & n\text{-point periodic Riesz minimizer, p. \pageref{p:periodic_minimizers}}\tabline 
            \lat_n', \ \rmin_n'                         & \tabbox{translated and rescaled $ \lat_n $ and $ \rmin_n $, \eqref{eq:irr_lat_scaled}-\eqref{eq:riesz_scaled}}\tabline
            n_m                                         & \tabbox{number of nodes in $ \ndv_m $} \tabline
            \Omega                                      & \tabbox{target distribution support} \tabline
            \Omega_\text{etopo},\ \Omega_\text{shell}   & \tabbox{the underlying sets in Sections \ref{subsec:atm_nodes} and \ref{subsec:sph_shell}} \tabline
            (\VarRadius,\,\VarAngle,\,\VarPolarAngle) & \tabbox{spherical coordinates, p. \pageref{p:spherical}}\tabline
            \rho(\vecx)                                 & \tabbox{objective radial density, pp. \pageref{p:radial_density},\;\pageref{p:radial_density_properties}}\tabline 
            \ndv_m,\ m=1,\ldots,\VoxelsPerSide^d        & \tabbox{cube-shaped voxel in $ \mathbb R^d $, p. \pageref{p:voxel}}\tabline
            \vecx; \ \vecx_i,\ i=1,\ldots,N             & \text{points in }\mathbb R^d;\text{ nodes of the configuration } \tabline
            \nd_{j(i,k)},\ k=1,\ldots,K            & \tabbox{the $k$-th nearest neighbor  of  $ \vecx_i $, \eqref{eq:sw-riesz}}\tabline 
            \bs{z}_m                                    & \tabbox{center of $\mathcal V_m$} \tabline
            \bottomrule
        \end{array}
    \]
    \caption{Symbolic notation employed throughout the paper.}
    \label{tab:notation}
\end{table}

\section{Choice of method}\label{sec:method} 
\subsection{RBF-FD approximations}\label{subsec:rbfs} In this section we shall outline the common practices involving RBFs, in order to motivate the requirements that have to be imposed on the node distribution used in the respective computations. For a more in-depth discussion see one of \cite{Fornberg2015b, Wendland2004, Buhmann2003, MR2357267}. A significant portion of the RBF approach hinges on the theory of positive definite functions.

Suppose we need to approximate a linear operator $ \mathfrak L $ acting on sufficiently smooth functions supported on $ \Omega $, given locally by their values at the nodes $ \bs x_k,\,k=1,\ldots,K $. More specifically, we need to compute the value $ \mathfrak L\psi(\bs x_0)$ for some fixed $ \bs x_0 \in \Omega $ and a variable function $ \psi $.   
A generalization of the standard \cite{Fornberg1998b} finite-difference (FD) approach consists in constructing weights $ w_k,\,k=1,\ldots,K  $, that recover the value of $ \mathfrak L $ at $ \bs x_0 $ in the form 
\begin{equation}
    \label{eq:weights}
    \mathfrak L\mathcal S(\bs x_0) = \sum_{k=1}^K w_k \mathcal S(\bs x_k), 
\end{equation}
for every $ \mathcal S $ from a convenient functional space; the $ \mathcal  S $ is then chosen to interpolate $ \psi $ at the given nodes $ \bs x_1,\ldots, \bs x_K $. In our case, $ \mathcal S $ is spanned by RBFs, so by analogy to the 1-dimensional case this method is called RBF-FD; there exists extensive literature covering different types of kernels and different applications \cite{Fornberg2015b, Flyer2013, FornbergLarssonFlyer, Bollig2012,Flyer2012,Flyer2016,Bayona2017,Flyer2016a}.  
Note that the node localization is required due to that (i) local stencils lead to sparse matrices, and are thus much more suitable for computations, (ii) in most applications, $ \mathfrak L $ is either an interpolation or a differential operator; both act locally, so it is natural to use local stencils.

For example, using the space of shifts of the Gaussian kernel $ \phi(r) = e^{-(\epsilon r)^2} $, one arrives at an RBF interpolant
\[
    \mathcal S(\bs x) = \sum_{k=1}^{K} c_k \phi(\|\bs x - \bs x_k\|).
\]
In order to express $ \mathfrak L\mathcal S(\bs x_0) $ as a functional of $ \mathcal S(\bs x_k), k=1,\ldots,K $, as in \eqref{eq:weights},
it suffices to do so for the functions $ \phi_k(\bs x) = \phi(\|\bs x - \bs x_k\|),\, k=1,\ldots,K  $. The weights $ \{w_k\} $ are then obtained as the solution to 
\[
    \begin{bmatrix} \phi(\|\vecx_1-\vecx_1\|)& \phi(\|\vecx_1-\vecx_2\|) &
        \ldots & \phi(\|\vecx_1-\vecx_K\|) \\ \phi(\|\vecx_2-\vecx_1\|)&
        \phi(\|\vecx_2-\vecx_2\|) & \ldots & \phi(\|\vecx_2-\vecx_K\|) \\
        \vdots& \vdots  &  & \vdots  \\ \phi(\|\vecx_K-\vecx_1\|)&
        \phi(\|\vecx_K-\vecx_2\|) & \ldots & \phi(\|\vecx_K-\vecx_K\|) \\
    \end{bmatrix}
    \begin{bmatrix}
        w_1\\
        w_2\\
        \vdots \\
        w_K
    \end{bmatrix} = 
    \begin{bmatrix}
        \mathfrak L\phi(\|\bs x_0 - \bs x_1\|)\\
        \mathfrak L\phi(\|\bs x_0 - \bs x_2\|)\\
            \vdots\\           
        \mathfrak L\phi(\|\bs x_0 - \bs x_K\|)\\
    \end{bmatrix}.
\]
Observe that in order to find the interpolant $ \mathcal S $ for $ \psi $, $ \{c_k\} $ are determined from the same system, with $ \mathfrak L $ taken to be the identity map.
The matrix on the left is denoted by $ \bs A $ in \eqref{eq:Amat}; it is degenerate whenever any two of $ \{\bs x_k\} $ coincide, and is ill-conditioned whenever any two are very close, due to the continuity of $ \phi $. This brings us to further considerations of how the stencil $ \{\bs x_k\} $ can be chosen. By the above, it is necessary that the nodes be (i) distinct and well-separated, and (ii) localized inside the domain $ \Omega $. For a quasi-uniform node set, $ K $ nearest neighbors of a fixed node satisfy both conditions.

Observe that for any strictly positive definite kernel $ \phi $, provided $ \{\bs x_k\} $ are all distinct, the interpolation matrix $ \bs A $ is always invertible (the Gaussian is an example of such kernel). 
To summarize, the above expression for weights $ \{w_k\} $ is the defining property of the RBF-FD methods with the Gaussian kernel.  

Taking the limit of the shape parameter $ \epsilon \to 0 $ can cause the interpolant $ s $ to diverge for other RBF kernels \cite{Fornberg2004b,Buhmann2010} that contain $ \epsilon $; this phenomenon however does not occur for the Gaussian $ \phi(r) = e^{-(\epsilon r)^2} $. The motivation for considering the ``increasingly flat'' limit $ \epsilon \to 0 $ is that the resulting RBFs can be used to obtain highly accurate solutions of elliptic problems and approximants of smooth data \cite{Fornberg2004b,Fornberg2013a}. We now conclude the discussion of the Gaussian kernel and turn to its novel alternative.

In the recent years, there have been noteworthy advances in RBF-FD using Polyharmonic Spline (PHS) kernels, $ \phi(r) = r^{2p-1} $ or $ \phi(r) = r^{2p}\log r,\, p\in \mathbb N $; it has has been shown \cite{Flyer2016,Bayona2017, Flyer2016a} that using PHS-based RBF-FD leads to improved accuracy, stability, and eliminates the Runge phenomenon at the boundary of the domain \cite{Bayona2017}, which is not the case in general \cite{Fornberg2007}. Another benefit from using the PHS kernel is that it does not contain the shape parameter $ \epsilon $. The analytical property underlying existence of the weights $ \{w_k\} $ for the PHS  kernels is that the PHS are conditionally strictly positive definite \cite{MR2357267, Micchelli1986} and thus need a slightly different treatment, which we shall now outline.

To ensure unisolvency (uniqueness of the weights and interpolants), we need to augment the $ \mathcal S $ with polynomial terms: it is selected from the space defined by
\[
    \mathcal S(\bs x) = \sum_{k=1}^K c_k \phi(\|\bs x - \bs x_k\|) + \sum_{i=1}^{\binom{l+d}{l}} b_i \pi_i(\bs x),
\]
with $ \{c_k\} $ satisfying the constraint 
\[
    \sum_{k=1}^K c_k \pi_i(\bs x_k) = 0, \quad i=1,2,\ldots, \binom{l+d}{l},
\]
where $ \phi $ is now the PHS kernel, and $ \binom{l+d}{l} $ is the dimension of the space of multivariate polynomials of degree up to $ l $ in $ \mathbb R^d $; accordingly, $ \pi_i $ varies over the monomial basis for such polynomials.  

The degree $ l $ has to satisfy $ l = p-1 $ and $ l = p $ for $ \phi(r) = r^{2p-1} $ and $ \phi(r) = r^{2p}\log r,$ respectively \cite[Chapter 8]{MR2357267}. For example, when $ \phi(r) = r^2\log r $ and $ d = 2 $, then $ l=1 $ and the weights corresponding to an operator $ \mathfrak L $ in $ \mathbb R^2 $ are determined by
\begin{equation}
    \label{eq:RBF-FD}
    \begin{bmatrix}
        & &         & |   & 1 & x_{11} & x_{12}       \\
        & \bs A &   & |   & \vdots & \vdots & \vdots  \\
        & &         & |   & 1 & x_{k1} & x_{k2}       \\
        - & - & -   & +   & - & - & - \\
        1 & 1 &  1  & | & & & \\
        x_{11} & \dots & x_{k1}  & | & &  0 & \\
        x_{12} & \dots & x_{k2}  & | & & & \\
    \end{bmatrix}
    \begin{bmatrix}
        w_1 \\
        \vdots\\
        w_k\\
        -\\
        w_{k+1}\\
        w_{k+2}\\
        w_{k+3}
    \end{bmatrix}
    =
    \begin{bmatrix}
        \mathfrak L\phi(\|\bs x_0 - \bs x_1\|)\\
            \vdots\\           
        \mathfrak L\phi(\|\bs x_0 - \bs x_K\|)\\
        -\\
       \mathfrak L\, 1\\
       \mathfrak Lx_1(x_{01})\\
       \mathfrak Lx_2(x_{02})
    \end{bmatrix}, 
\end{equation}
where the matrix $ \bs A $ is the same as in \eqref{eq:Amat} with the PHS kernel; $ x_{kj} $ is the $ j $-th coordinate of $ \bs x_k,\, k=1,\ldots,K;\, j=1,2 $, and similarly for $ \bs x_0 = (x_{01}, x_{02})^\text{tr} $; $ \mathfrak L1,\, \mathfrak L x_1,\, \mathfrak L x_2 $ denote images of the constant and coordinate functions under $ \mathfrak L $, respectively. Here, as before, the interpolation case is obtained by taking $ \mathfrak L $ equal to the identity operator; compare the constraints on $ \{c_k\} $ above. The generalization to larger values of $ l $ and higher dimensions follows along the same lines, with higher-degree monomials instead of linear terms \cite[Chapters 8, 11.1]{Wendland2004}.  

For this and the other commonly used kernels, non-degeneracy follows from a strengthened form of Micchelli's theorem \cite{Powell2004}, see also \cite{MR2357267, Wendland2004, MR2474372,Micchelli1986}; namely, the matrix in the LHS of the previous equation is non-degenerate for any unisolvent $ \bs x_1,\ldots,\bs x_K $. The remaining part of the discussion for the Gaussian kernel above is further applicable without any modifications. It should be noted that the optimal choice of the degree of PHS that needs to be used, does in general depend on the particular problem under consideration.

\subsection{Riesz energy}\label{subsec:riesz} To generate nodes both devoid of lattice alignment 
and having near-optimal local separation, we shall minimize a functional on the space of discrete subsets of $ \Omega $. Equivalently, one can think of the corresponding gradient flow moving the starting configuration to a suitable position. The desired properties of the minimizing configuration will then follow from the strong repulsion imposed by the functional.

First, for a fixed $ s > d $ we introduce \textit{Riesz $ s $-energy}, a functional $ E_s:\Omega^N\mapsto \left(0,\infty \right) $ such that for a collection of vectors $ \vecx_1,\ldots,\vecx_N $ in $ \Omega $,
\begin{equation}\label{eq:s-riesz}
E_s(\vecx_1,\ldots\vecx_N) = \sum_{i\neq j} \frac{1}{ \|\vecx_i -
\vecx_j\|^s}.  
\end{equation} There exists extensive literature dedicated to the collections minimizing this and derived functionals for $ s\geq d $, their asymptotics and limiting measures, see for example \cite{Hardin2005, Borodachov2008, Brauchart2012b}.
It turns out that in the case $ s\geq d $ the distribution of minimizers of $ E_s $ coincides with the normalized Hausdorff measure on $ \Omega $; practically this  means that the minimizers are uniform in the volumetric sense on $ \Omega $, that is, the number of nodes per unit volume is close to constant. In order to produce non-uniform nodes, we shall further add  multiplicative weight to \eqref{eq:s-riesz}; this modification of \eqref{eq:s-riesz} was first studied in \cite{Borodachov2008}. The \textit{weighted Riesz
$ s $-energy} with kernel $ \kappa: \Omega\times\Omega \to (0,\infty) $ is the functional $ E_s^\kappa:\Omega^N\mapsto \left(0,\infty \right)  $ defined by
\begin{equation*}
    E_s^\kappa(\vecx_1,\ldots\vecx_N) = \sum_{i\neq j}
\frac{\kappa(\vecx_i , \vecx_j)}{ \|\vecx_i - \vecx_j\|^s}.
\end{equation*} 
It has been shown \cite{Borodachov2008}, that the counting measures of the minimizers of the weighted energy converge to the probability measure with volumetric density proportional to $ \kappa(\vecx, \vecx)^{-d/s} d\mathcal{H}_d(\vecx) $, with $ \mathcal{H}_d $ denoting the $ d $-dimensional Hausdorff measure. More precisely, for any $ \mathcal B \subset \Omega $ with boundary of zero measure there holds
\[
    \frac1N\sum_{i=1}^N \chi(\bs x_i; \mathcal B) \longrightarrow \frac1{Z(\kappa, \Omega)}\int_\mathcal B\kappa^{-d/s}(\bs x, \bs x)d\mathcal H_d(\bs x), \quad N\to \infty,
\]
where $ \chi(\cdot;\mathcal B) $ is the indicator function of $ \mathcal B $ and $ Z(\kappa, \Omega) $ is the normalization constant. Of course, for most applications the set $ \mathcal B $ will have zero-measure or even differentiable boundary.

It is worth noting that the previous equation shows that the distribution of nodes produced by minimizing $ E_s^\kappa $ depends only on the diagonal values of $ \kappa $ for large enough $ N $. Indeed, this has been explored in \cite{Borodachov2014}, where it is shown that omitting interactions of points at least $ r_N $ apart in $ E_s^\kappa $ and minimizing the resulting expression leads to the same distribution when $ N\to \infty $; the sequence $ r_N,\, N\geq 2 $ here satisfies $ r_N N^{1/d} \to \infty,\, N\to \infty $. Following \cite{Borodachov2014}, a weight $ \kappa(\bs x,\bs y) $ vanishing whenever $ \|\bs x - \bs y\| > r_N $ is said to be \textit{truncated}.

Configurations that minimize $ E_s^\kappa $ over $ \Omega^N $ for a compact $ \Omega $ are well-separated, that is, the quantity  $ \Delta(\{\vecx_1,\ldots\vecx_N\}) = \min_i\Delta(\vecx_i) $ satisfies 
\begin{equation}\label{eq:separation}
    \Delta(\{\vecx_1,\ldots\vecx_N\}) \geq C N^{-1/d}, \quad N\geq 2.
\end{equation} 
We shall assume that the terms in $ E_s^\kappa $ for which $ \vecx_j $ is not among the $ K $ nearest nodes to $ \vecx_i $ are zero, a condition equivalent to truncating $ \kappa $, provided the nodes are well-separated. Under this assumption, the expression for $ E_s^\kappa $ can be rewritten as
\begin{equation}
    \label{eq:sw-riesz}
    E_s^\kappa(\vecx_1,\ldots\vecx_N) = \sum_{i=1}^N \sum_{k=1}^K \frac{\kappa\big(\vecx_i , \nd_{j(i,k)}\big)}{ \|\vecx_i - \nd_{j(i,k)}\|^s},
\end{equation}
where nodes $ \nd_{j(i,k)}, $ $ k=1,\ldots,K $, are the $ K $ nearest neighbors of $ \vecx $. That minimizers of \eqref{eq:sw-riesz} are well-separated can be shown by the standard argument from one of \cite{Hardin2005,Borodachov2008,Borodachov2014}. This further implies that they are quasi-uniform, which is the key property for us in view of the discussion in Section~\ref{subsec:rbfs}. As the form \eqref{eq:sw-riesz} makes clear, for the truncated kernel $ \kappa $ the $ E_s^\kappa $ can be
computed in $ O(NK) $ operations, unlike the $ O(N^2) $
operations required to evaluate the functional for a non-vanishing $ \kappa $. This, and the fact that \eqref{eq:sw-riesz} requires constant size memory for storage makes this form of $ E_s^\kappa $ the most useful for our purposes.

The value of the exponent $ s $ is
chosen so that $ s\geq d $ to ensure that the energy functional is sufficiently
repulsive; it is known from the classical potential theory that for $ s<d $ the
minimal energy configurations are not necessarily uniform, and their local
structure depends on the shape of the domain \cite{Landkof1972}. Property
\eqref{eq:separation} holds for any $ s>d $, when minimizing the energy over any
fixed compact set $ \Omega \subset \mathbb{R}^d $. While setting $ \kappa(\bs x, \bs y) = f(\|\bs x - \bs y \|) $ for an arbitrary positive definite radial function $ f(r) $ that grows fast enough when $ r\to 0 $ would produce similar results, we chose the Riesz energy because the properties of its minimizers are well understood.

Note that simply looking for minimizers of $ E_s^\kappa $ does not provide node sets satisfying $ \Delta(\nd) \asymp \rho(\nd) $ for every $ \bs x $; in fact, boundary nodes of such minimizers will often have smaller separations than desired. Since in such cases the boundary has a lower Hausdorff dimension, it does not influence the volumetric density, which agrees with the results above. With this motivation in mind, we are ready to introduce the second component of our method.

\subsection{Quasi-Monte Carlo methods}\label{subsec:qmc}
To facilitate convergence of whichever optimization algorithm is used to find minimizers of \eqref{eq:sw-riesz}, we can initialize it with a configuration that approximates the limiting measure. One has to rule out Monte Carlo methods due
to the separation requirement: 
random points exhibit clustering \cite{Brauchart2016}, which makes deterministic post-processing, in particular by energy minimization, costly. Similarly, mitigating the clustering by purely probabilistic approaches, as for example \textit{thinning} discussed in \cite{Link2012}, or by sampling from a random process with repulsive properties \cite{Alishahi2015a,Beltran2016a}, does not generally yield satisfying results, since the separation can only be guaranteed on average. Instead, we turn to the
quasi-Monte Carlo (Q-MC) approach. As has been pointed out at the end of Section~\ref{subsec:rbfs}, in order to ensure convergence of RBF-FD interpolants, the underlying node set must be (locally) unisolvent; for our purposes this just means that the nodes are in a generic position with respect to each other. The latter is clearly not the case for lattice-like Q-MC configurations, which explains why we resort to energy minimization. On the other hand, we choose not to use other popular Q-MC sequences, such as Halton nodes, as they do not necessarily lead to the best conditioned systems, see \cite[Figure 5.1]{Fornberg2013} and Figure~\ref{fig:condition_numbers} below, and are harder to handle when the distribution support $ \Omega $ has complex geometry.

The key element of our construction lies in distributing the node set in a
deterministic way so as to guarantee low discrepancy between the desired and
the obtained radial densities. This is achieved by a Q-MC analog of the
\textit{stratification} of the Monte Carlo method \cite{Caflisch1998}: nodes are
distributed with piecewise constant (radial) density that approximates the
desired one. We consider two different Q-MC sequences to draw from with
(near-)constant radial density: irrational lattices and periodic Riesz
minimizers. After dividing the set $ \Omega $ into cube-shaped \label{p:voxel}\textit{voxels},
each voxel is filled with nodes obtained in one of the two ways, appropriately
scaled, then the weighted $ s $-energy \eqref{eq:sw-riesz} of the whole node set is minimized.  Although we discuss the radial density in the present paper,
an argument for the volumetric density can be produced along the same lines.

Yet another reason to make use of a  Q-MC sequence is to avoid
recursive  data structures, which in some cases can be detrimental to the
overall performance. Even though such structures have seen significant developments over the years, both
dynamic update and parallelization for them remain challenging, 
\cite{Procopiuc_2003, Zhou2008}. The approach of the present paper should therefore be understood as almost opposite to the well-known ``quadtree" algorithm \cite{FreyGeorge2010},  that indeed has been used for meshless node generation \cite{Varma2004}. Namely, as outlined above, our algorithm places nodes en masse inside the voxels to produce a rough approximation of the target distribution, and subsequently adjusts them by a gradient flow, which is straightforward to parallelize. Although this does involve the computation of the nearest neighbors in \eqref{eq:sw-riesz}, which in practice will be done by constructing a k-d tree, by initializing the node configuration with a stratified Q-MC sequence we ensure the indices $ j(i,k) $ in \eqref{eq:sw-riesz} will not undergo significant changes during the energy minimization stage, so the k-d tree will not require intensive updates.

An \textit{irrational lattice} (IL) is defined as a discrete subset of the
$ d $-dimensional unit cube $ [0,1]^d $ 
\begin{equation}\label{eq:irr_lat} 
    \lat_n
    = \bigg\{\bigg(\Fr(\OffsetLattice+{i}/n), \Fr(i\alpha_1), \{ i\alpha_2   \},
    \ldots, \Fr(i\alpha_{d-1})\bigg) \ \bs{:}\ i=1,\ldots,n\bigg\}, 
\end{equation}
where $\Fr(x) = \mathrm{mod}(x,1) = x - \lfloor x \rfloor $ denotes the fractional part of
$x$, $ \OffsetLattice>0 $ is fixed, and $\alpha_1, \alpha_2, \ldots,
\alpha_{d-1}$ are irrational numbers, linearly independent over the
rationals. This terminology seems to be accepted in the low-discrepancy
community \cite{Bilyk2013}, while closely related objects, when used
for Q-MC purposes, are known as \textit{Korobov/lattice point sets}
\cite{Lemieux2009}. 

The motivation for using an IL in this context is due to the
existing results on the discrepancy of ILs. It is known for example, that the
two-dimensional ILs have the optimal order of $ L^2 $ discrepancy,
\cite{Bilyk2012, Bilyk2013}. Furthermore, in all dimensions ILs are
uniformly distributed \cite[Chapter 1.6]{Kuipers2006}, that is, the fraction of
lattice points inside any rectangular box with faces parallel to the coordinate
planes converges to its volume. The simple linear structure of ILs makes them
especially attractive for SIMD-parallelization.

Another Q-MC sequence that has proven to suit our purposes consists of
\label{p:periodic_minimizers}\textit{periodic Riesz minimizers} on the unit flat torus, that is, $ n $-point collections $ \mathcal M_n = \{\bs{x}_1,\ldots,\bs{x}_n \} $ that minimize
 \eqref{eq:s-riesz} on $ ([0,1]^d)^n $ with the Euclidean distance $ \|\cdot\| $ replaced by the
periodic metric 
\begin{equation}\label{eq:pdist}
    \|\bs{x} - \bs{y}\|_\sim^2 = \Pi(x_1-y_1)  + \Pi(x_2-y_2)  +
    \Pi(x_3-y_3),  
\end{equation}
where $ \Pi(x) = \min(x^2,  (1-x)^2 ), \ 0\leq x\leq 1 $. It
follows from \cite{Hardin2005} that such configurations have optimal separation
and asymptotically uniform volumetric density. It follows from the numerics
also, that in this case the nearest neighbor distances vary very little from
node to node; this and that minimizing configurations do not suffer from the
lattice-like alignment, makes their rescaled copies good candidates for the
stratification.

The number of nodes in individual voxels is defined by the function $ \rho $,
so the resulting collection has piecewise constant density; refining the voxel
partition leads to an improved piecewise approximation of the desired (e.g.,
smooth) density. In practice, the dependence between the number of nodes
contained in the unit cube, and average/minimal nearest neighbor distance is
tabulated in advance, and then inverted during the construction of the node
set.

\begin{figure}[t]
    \centering
    \includegraphics[width=\smwdth\textwidth]{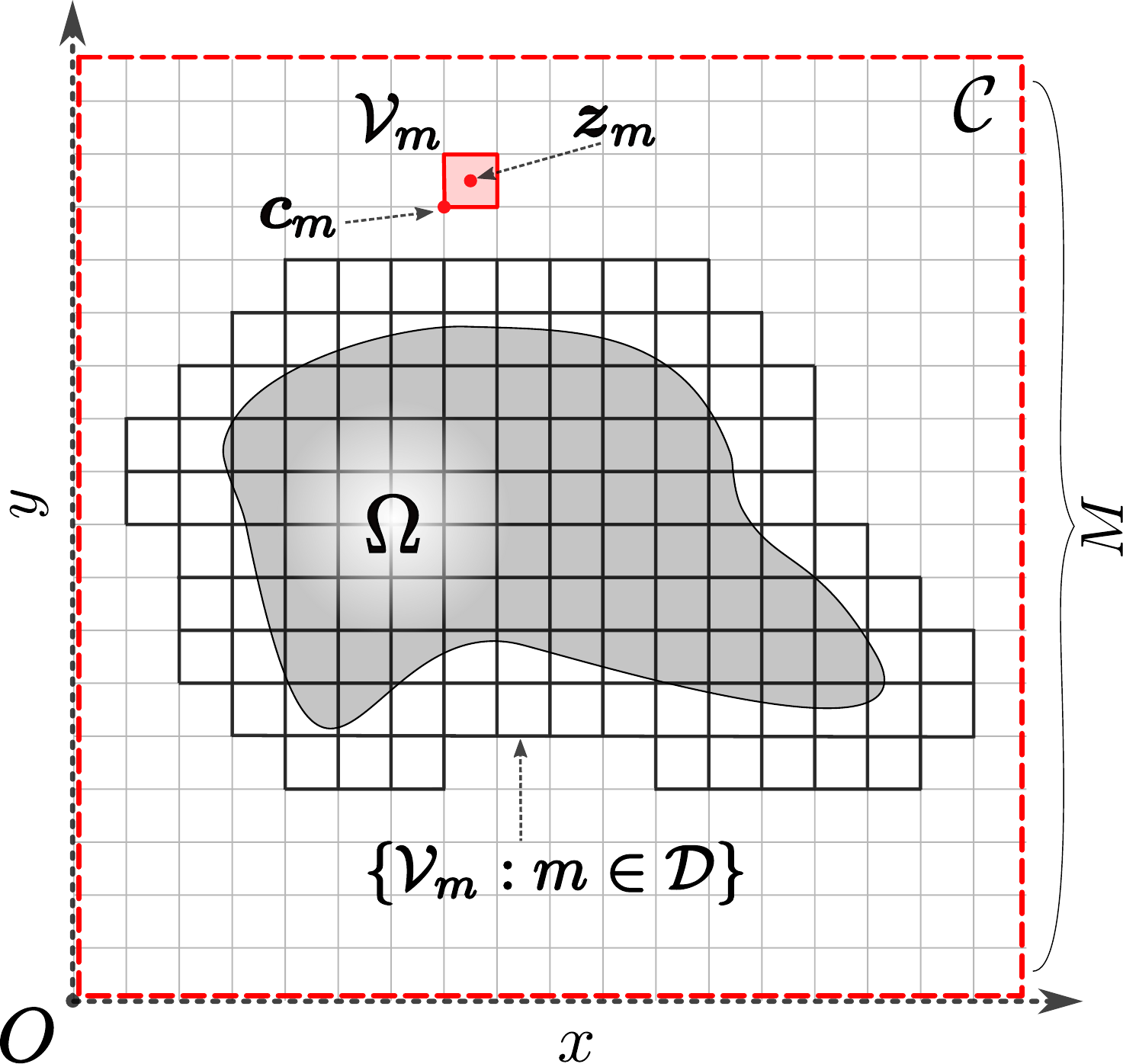}
    \caption{An illustration of some of the symbolic notation used in the algorithm below, in the case $ d=2 $.}
    \label{fig:notation}
\end{figure}
\section{The algorithm}\label{sec:algo}
The interested reader will find a Matlab codebase implementing the algorithm described here, as well as the sources for all the figures contained in the present paper,  at \cite{VMFF17matlab}.
\subsection{Formulation}\label{subsec:algo}
If the nodes must be restricted to a certain compact
set $ \Omega $, for example, support of a given indicator, we will refer to the set as
\textit{density support}, and to the indicator function as \textit{point
inclusion function}. We may further assume that $ \Omega $ is
contained in the $ d $-dimensional unit cube $ \mathcal{C} =[0,1]^d $ (see Figure~\ref{fig:notation} for some of the notation involved); the case of an arbitrary compact
set then follows by choosing a suitable enclosing cube and applying scaling
and translation. Suppose the radial density is prescribed by a Lipschitz-1
function; i.e., \label{p:radial_density_properties}$ |\rho(\bs{x})-\rho(\bs{y})|\leq \|\bs{x}-\bs{y}\| $. The reason for this assumption is the respective property of $ \Delta(\nd) $, and is explained in further detail in the following section. Recall that we use an exponent
$s>d$. We summarize the discussion in Section~\ref{sec:method} into the
following algorithm for generating nodes with radial density $ \rho $:  

\vspace*{\baselineskip}
\noindent \textbf{Initial node layout.}
\hrule
\begin{algo}
\item \label{step:prepare} Choose one of the two Q-MC sequences
    described in Section~\ref{sec:method}, $ \{\lat_n : n\geq 1\} $ or $ \{\rmin_n : n~\geq~1\} $,
    draw configurations with up to $
    n_{\max} $ nodes from it, and determine the average nearest neighbor distance
    for its periodization by the integer lattice, denoted by \label{p:delta_n}$  \bar\Delta_n $ for $ n $ nodes. Let $ \lambda: (0,\infty) \to \{0,1,2,\ldots, n_{\max}\} $ be the interpolated
    inverse to $ \bar\Delta_n: \{1,2,\ldots, n_{\max}\} \to (0,1] $.\footnote{Note that
        both ILs and the minimizers can have the nearest neighbor distance of at most $
    1 $, due to periodicity. We therefore take $ \lambda(x)=0 $ whenever $ x>1 $.} \item
    \label{step:subcubes}  Partition $\mathcal{C}$ into $\VoxelsPerSide^d$ equal cube-shaped voxels $\ndv_m,$  $m=1,\ldots,\VoxelsPerSide^d$ of side length $1/\VoxelsPerSide$, with
    faces parallel to the coordinate planes.  Let \label{p:D_nonempty}$ \{\ndv_m : m\in\mathcal{D}\} $
    be the subset for which at least one of the adjacent (i.e., sharing a face)
    voxels has a vertex inside $ \Omega $.

\item \label{step:densefill}  Let $\bar\rho_m$ be the average value of $\rho$
    at the $ 2^d $ vertices  of a voxel $ \ndv_m,\ m\in\mathcal{D} $. Place
    inside $ \ndv_m $ a scaled and translated version of the $ n_m $-point IL
    \eqref{eq:irr_lat}, or of the $ n_m $-point periodic Riesz minimizer, using $
    n_m $ defined by \[  n_m = \lambda(\bar{\rho}_m  \VoxelsPerSide ).  \]
    Repeat for each $ m\in\mathcal{D}  $.
\end{algo}
\noindent \textbf{Saturation and cleanup.}
\hrule
\begin{algo2}
\item\label{step:saturate}
    Let $ \mathcal E \subset \mathcal D $ be the set of $ m $ for which voxels $ \{\ndv_m : m\in\mathcal{E}\} $ contain no nodes and the centers $ \{\bs z_m : m\in\mathcal{E}\} $ satisfy $ \Delta(\bs{z}_m) > \rho(\bs z_m) $. Sort $ \mathcal E $ by the increasing values of $ \rho(\bs z_m) $. Repeat until $ \mathcal E $ is empty: for every $ m \in \mathcal E $  place a node in $ \bs{z}_m $; recompute $ \mathcal E $.

\item\label{step:cleanup} For all nonempty voxels, remove nodes outside $ \Omega $.
\end{algo2}
\noindent \textbf{Repel-type iterations, boundary detection.}
\hrule
\begin{algo3}
\item\label{step:repel} Perform $ T $ iterations of the partial gradient
    descent on the weighted $ s $-energy functional \eqref{eq:sw-riesz} with $
    \kappa(\vecx,\vecx) = c\rho(\vecx)^s $, using the $ K $ nearest neighbors of
    each node:
    Let the initial configuration be the $ 0 $-th iteration, $
    \nd_i^{(0)}= \nd_i,$  $ i=1,\ldots,N  $, with $
    N $ denoting the total number of nodes distributed.  On the $t^{th}$
    iteration, $ 1\leq t\leq T $, given a node $\nd_i^{(t)}$ with $
    K $ nearest neighbors $ \nd_{j(i,k)}^{(t)}, $ $ k=1,\ldots,K $, form
    the weighted vector sum 
    \[
        \bs{g}_i^{(t)} =
        s\,\rho\left(\nd_i^{(t)}\right)^s
        \sum_{k=1}^{K}\frac{\nd_i^{(t)}-\nd_{j(i,k)}^{(t)}
        }{\|\nd_i^{(t)}-\nd_{j(i,k)}^{(t)} \|^{s+2}}, \quad1\leq i \leq N,
    \]
    and the new
    node position can now be expressed as 
    \begin{equation}\label{eq:node_update}
        \nd_i^{(t+1)} = \begin{cases} \displaystyle \nd_i^{(t)} +
            \frac{\Delta(\nd_i^{(t)})}{t + \OffsetGrad}  \frac{ \bs{g}_i^{(t)}
            }{\|\bs{g}_i^{(t)}\|} &\text{if this sum is inside }\Omega;
        \\ \nd_i^{(t)}, &\text{otherwise,}	\end{cases} \quad1\leq i \leq N.
    \end{equation}
    where $\OffsetGrad$ is a fixed offset chosen to control the step
    size between $\nd_i^{(t)}$ and $\nd_i^{(t+1)}$. If a ``pullback'' function
    is provided from a neighborhood of $ \Omega $ to its boundary, the
    condition of $\nd_i^{(t+1)}$ being inside $ \Omega $ is
    replaced with applying the pullback; furthermore, if the radial density has
    an easily computable gradient, or is changing rapidly, an additional term
    must be included in \eqref{eq:node_update} (see discussion below). \\Update
    the neighbor indices $ j(i,k) $ after every few iterations.

\item\label{step:surface_nodes} If no boundary node set/pullback function is
    prescribed, define the boundary nodes as follows. Evaluate the point
    inclusion function for $ \nd_i \pm \Delta(\nd_i) \bs{e}_l, \ l=1,\ldots,d,
    \ i=1,\ldots,N $, where $\bs{e}_l $ is the $l$-th basis vector. If at least
    one such point lies outside  $ \Omega $, the $\nd_i$ is considered to be a boundary node.  
\end{algo3} 
\subsection{Discussion}\label{subsec:algo_discussed}
Our assumption of $ \rho $
being Lipschitz-$ 1 $ is natural, since $ \Delta(\cdot) $ is always Lipschitz-$
1 $, if viewed as a function of position. To see this, consider any two nodes $
\bs{x}, \bs{y} $, and let $ \bs{x'}, \bs{y'} $ be their nearest neighbors,
respectively, so that $ \|\bs{x} - \bs{x'}\| = \Delta(\bs{x}) $ and $ \|\bs{y}
- \bs{y'}\| = \Delta(\bs{y}) $.  It follows, $ \Delta(\bs{x})\leq \|\bs{x} -
\bs{y'}\| \leq \| \bs{x}-\bs{y}\| +  \Delta(\bs{y})$, which by symmetry implies
$ \Delta $ is Lipschitz-1.

\subsubsection*{Initial node layout.} \label{sub:initial_node_layout} 
In the parts of the density support with constant $ \rho $, the nodes will
locally look like a periodization of the initial Q-MC sequence; hence the
average neighbor distance in \ref{step:prepare} is tabulated for the
periodized version. Observe that there is some freedom in the notion of voxel adjacency used to define $ \{\ndv_m : m\in\mathcal{D}\} $ in \ref{step:subcubes}; for example, in the case of a highly non-convex domain $ \Omega $, it might be reasonable to denote the $ 3^d-1 $ voxels sharing a vertex with a given $ \ndv_m $ as adjacent to it, rather than only the $ 2d $ voxels that it has a common face with. This would then ensure that no part of $ \Omega $ will be omitted in the node allocation; imagine a long and thin peninsula in $ \Omega $ containing no corners of $\ndv_m,$  $m=1,\ldots,\VoxelsPerSide^d$. We have found however, that the subsequent repel iterations will guarantee that such a peninsula is adequately filled with nodes even when using only the face-adjacent voxels.

If the IL sequence is chosen
in \ref{step:prepare}, the $ n_m $-node set placed in voxel $\ndv_m$ at \ref{step:densefill} is an
adjusted version of \eqref{eq:irr_lat} as follows. Let for every $\ndv_m$ the corner with the smallest absolute value be $\bs{c}_m$; the points $ \bs{c}_m $ are then vertices of a lattice. Before scaling and translating $ \lat_{n_m} $, apply a random permutation  to the coordinates of each node in it, so
that to remove long-range lattice structure from the distribution; we will denote such an operation by $\sigma$. Then the IL  in voxel $\ndv_m$ becomes
\begin{equation}\label{eq:irr_lat_scaled} 
    \lat_{n_m}'  =  \bs{c}_m + \frac{\ScalingFactor}{\VoxelsPerSide} \sigma(\lat_{n_m})  + \frac{\ShiftVector}{\VoxelsPerSide}, 
\end{equation}
where
\[
\begin{aligned}
    \ScalingFactor      =  { 1 - c_d\left({n_\text{max}}\right)^{-1/d} },
    \qquad
    \ShiftVector =  \frac{1-\ScalingFactor}{2}\cdot (1,1,\ldots,1)^\text{tr},
\end{aligned}
\]
with $ c_d $ depending only on the dimension. The
quantities $\ScalingFactor$ and $\ShiftVector$ ensure that the lattice points in $\lat_m$ are
inset into the voxel by about half the separation distance, avoiding
poorly separated points along the voxel interfaces.

When the periodic minimizer sequence is selected in \ref{step:prepare}, the inset is defined in a similar way, but the permutation is just an identity, $ \sigma \equiv \text{id} $, as the minimizers don't have the lattice structure. Likewise, the scaling factor and translation are
\[
\begin{aligned}
    \ScalingFactor_m      =  { 1 - c_d\left({n_m}\right)^{-1/d} },
    \qquad
    \ShiftVector_m =  \frac{1-\ScalingFactor_m}{2}\cdot (1,1,\ldots,1)^\text{tr}.
\end{aligned}
\]
The analog of \eqref{eq:irr_lat_scaled} thus takes the form
\begin{equation}\label{eq:riesz_scaled} 
    \rmin_{n_m}'  =  \bs{c}_m +
    \frac{\ScalingFactor_m}{\VoxelsPerSide}  \rmin_{n_m}  + \frac{\ShiftVector_m}{\VoxelsPerSide}.
\end{equation}
As one would expect, the average separation for the sequence $ \{\rmin_n\} $ is larger than that of $ \{\lat_n\} $ for the respective values of $ n $. While the inset for the latter is necessary to account for the node proximity after periodization, for the former it serves to mitigate the effects of interfacing voxels containing different number of nodes. This is further discussed in the Appendix.

\subsubsection*{Saturation and cleanup.}
\label{sub:saturation_and_cleanup} 
Observe that after \ref{step:densefill}, voxels in $ \{\ndv_m : m\in\mathcal{D}\} $ satisfying $ \bar{\rho}_m  \VoxelsPerSide > 1 $ do not contain any nodes. The goal of \ref{step:saturate} is therefore to remove any redundant sparsity that may be present whenever the radial density $ \rho $ is larger than $ 1/M $, as in this case the  function $ \lambda $ in \ref{step:densefill} is set to zero. More careful geometric considerations would lead one to set $  \lambda(x) > 0 $ when $ 0 < x <  \sqrt d/M  $, the length of a voxel diagonal, and thus make $ \lambda $ dependent on the dimension; on the other hand, using the interval $ 0 < x < 1/M $ as we did appears to suffice due to correction of density in \ref{step:repel}. 

Note that in practice, when recomputing $ \mathcal{E} $ in \ref{step:saturate}, to verify $ \Delta(\bs z_{m_0}) > \rho(\bs z_{m_0}) $ for a fixed $ m_0 \in\mathcal{E} $ it is enough to check $ \| \bs z_{m_0} - \bs z_m \| > \rho(\bs z_{m_0}) $ for the previously selected $ \bs z_m $ with $ \rho(\bs z_m) < \rho(\bs z_{m_0}) $.
Indeed, let  $ \bs z_{m_0} $ be the center of $ \ndv_{m_0} $. Then, by the definition of $ \lambda $ in \ref{step:densefill}, the radial density $ \rho(\bs{z}_{m_0}) = (1+D)/M $ for some $ D>0  $, so the Lipschitz-1 property implies, for any $ \vecx $ such that $ \rho(\vecx) \leq 1/M $ there holds $ \| \bs{z}_{m_0} - \vecx  \| \geq |\rho(\bs z_{m_0}) - \rho(\vecx)| \geq D/M $. This ensures that distances from $ \bs z_{m_0} $ to the nodes produced on \ref{step:densefill} satisfy
\[
    \| \bs z_{m_0} - \bs x \| \geq \frac D{1+D}\rho(\bs z_{m_0}). 
\]
This shows, when $ D \gg 1 $, not checking the inequality $ \| \bs z_{m_0} - \bs x \| > \rho(\bs z_{m_0}) $ leads to at most a bounded factor error. On the other hand, for $ D \ll 1 $ distances from $ \bs z_{m_0} $ to the nodes from \ref{step:densefill}  are also controlled: it follows from \eqref{eq:irr_lat_scaled}--\eqref{eq:riesz_scaled} that for $ n_m $ small, nodes in the voxel $ \ndv_m $ have larger inset (depending on $ c_d $).  
This analysis is certainly not rigorous; however, applying the partial gradient descent in \ref{step:repel}, we are able to ensure that the ratio $ \rho/\Delta $ is close to~1, as desired.  

Observe that in \ref{step:densefill} the nodes are only placed in $ \ndv_m $'s for which either of the adjacent voxels has corners inside the density support $ \Omega $, so removing nodes outside $ \Omega $ in \ref{step:cleanup} does not lead to much overhead. Furthermore, since the density is evaluated at the corners only and not at individual nodes, the total number of evaluations may be significantly reduced, which is especially useful when  $ \rho $ is computationally expensive. It is essential here that due to the Lipschitz-1 property, $ \rho $ is well estimated by its values at the corners $ \bs c_m $; specifically, $ |\rho(\bs x) - \rho(\bs c_m)|\leq \sqrt d/2M $ with $ \bs c_m $ the nearest voxel corner to $ \bs x $.

\subsubsection*{Repel-type iterations, boundary detection.}
\label{sub:repel_type_iterations_boundary_node_detection}

The equality $ \kappa(\vecx,\vecx) = c\rho(\vecx)^s $ can be justified by observing that each node $ \nd $ of the target distribution must be contained in a ball of radius $ \rho(\nd) $, not containing any of the other nodes, hence, the volumetric density must be inverse proportional to $ \rho(\vecx)^d $. On the other hand, minimizers of \eqref{eq:sw-riesz} converge to the distribution with volumetric density $ \kappa(\vecx,\vecx)^{-d/s} $; hence $ \kappa(\vecx,\vecx)^{d/s} = c \rho(\vecx)^d  $.

The vector sum in \ref{step:repel} is the partial
$ \bs x $-gradient of  the weighted  Riesz $ s $-energy \eqref{eq:sw-riesz} in the sense that a single summand of  \eqref{eq:sw-riesz} is  $
e(\bs{x},\bs{y}) = \kappa(\bs{x}, \bs{y})\|\bs{x}-\bs{y}\|^{-s} $, and thus its complete $ \bs x $-gradient is equal to 
\[
    \nabla_{\bs{x}}\, e(\bs{x},\bs{y})  = -s\, \kappa (\bs{x},\bs{y}) (\bs{x}-\bs{y})
    \|\bs{x}-\bs{y}\|^{-s-2} + \nabla_{\bs{x}} \kappa (\bs{x},\bs{y})
\|\bs{x}-\bs{y}\|^{-s}.  
\]
For our purposes, the $ \bs{y} $ here is one of the $ K $ nearest nodes to $ \bs{x} $, and, since due to the Q-MC initialization there will be few isolated nodes, and since
off-diagonal values of $ \kappa(\bs{x}, \bs{y}) $ do not influence the limiting
distribution (for details see \cite{Borodachov2014}), we assume $ \kappa(\bs{x},
\bs{y}) \approx \kappa(\bs{x}, \bs{x}) $ to rewrite the previous equation as 
\begin{equation}\label{eq:fullgradient}
    \nabla_{\bs{x}}\, e(\bs{x},\bs{y})  =
    -s\,\rho\left(\boldsymbol{x}\right)^s (\bs{x}-\bs{y})
    \|\bs{x}-\bs{y}\|^{-s-2} + s\nabla_{\bs{x}} \rho(\bs{x}) \,
    \rho(\bs{x})^{s-1} \|\bs{x}-\bs{y}\|^{-s}.  
\end{equation}
As has been pointed out at the beginning of Section~\ref{subsec:algo_discussed}, in order to be meaningful as a radial density, the function $ \rho $ must be Lipschitz-1. Then by the Rademacher's theorem, $ \nabla_{\bs{x}} \rho $ exists almost everywhere; this validates the use of it in \eqref{eq:fullgradient} as well as the approximation $ \kappa(\bs{x},
\bs{y}) \approx \kappa(\bs{x}, \bs{x}) $.
The ratio of the second term to the first one in \eqref{eq:fullgradient} is bounded by $ \nabla_{\bs{x}}
\rho(\bs{x})\|\bs{x}-\bs{y}\|/\rho(\bs{x})  $ and, provided that the
distances from $ \bs{x} $ to its nearest neighbors are close to the
value of $ \rho(\bs{x}) $, is at most $ c\nabla_{\bs{x}} \rho(\bs{x}) $
for a constant $ c $. This condition is satisfied because the chosen
Q-MC sequences have very regular local structure. In practice, the node
distance is small on the scale of the support and varies slowly, so the
second term will have negligible impact on the direction of the
gradient after normalization; besides, precise gradient computation may
prove costly. For these reasons we omit the second term in equation
\eqref{eq:node_update}. If it is necessary to deal with a fast-changing
radial density, a trade-off between the computational costs and the
resulting distribution properties must be sought.

It doesn't matter which minimization method is applied to the weighted $ s
$-energy, rather the gradient descent is chosen due to its simplicity. Note,
the second case in \eqref{eq:node_update}, leading to shrinking of the line
stepping distance, can be thought of as a simplistic backtracking line search;
it turns out to be sufficient for our purposes.  Furthermore, applying a more
involved line search may significantly degrade performance for complicated or
nonsmooth domains.

The number of nearest neighbors $ K $ in \eqref{eq:sw-riesz}, and the number of iterations $ T $ in \ref{step:repel} can be adjusted to
achieve a trade-off between execution speed/memory consumption/local
separation. In our experiments,\footnote{The Matlab code we provide performs naive
    autotuning of $ K $ and $ T $, using the total number of nodes to be
    placed. Although sufficient for demonstration purposes, there is room for
improvement.} even relatively small values of $ K $ and $ T $  produce good
results: we used $ {K\approx T\approx 30} $ for 1.36 million nodes with constant density in Section~\ref{subsec:atm_nodes}, and $ K= 30,\, T= 200 $ for 0.58 million and 0.36 million nodes with variable densities in Sections~\ref{subsec:cloud} and \ref{subsec:sph_shell}, respectively.

\section{Sample applications}
\subsection{Atmospheric node distribution using surface
data}\label{subsec:atm_nodes}

We use the geodata \cite{ETOPO1} from the
collection of global relief
datasets produced by NOAA (National Oceanic and Atmospheric Administration), which
contains a 1 arc-minute resolution model. We generate a sample configuration consisting of 
1,356,566 nodes distributed uniformly inside an atmospheric-type shell $ \Omega_\text{etopo} $: the
outer boundary of $ \Omega_\text{etopo} $ is spherical, the inner one is an interpolation of
the relief from ETOPO1 data, exaggerated by a factor of $ 100 $. The scale is
chosen so that the average Earth radius, assumed to be 6,371,220 meters, has
unit length; the radius of the outer boundary is set to $ 1.1 $, which corresponds to the height of 6,371 meters above the average radius, given the exaggeration factor.

\begin{figure} \centering 
    \begin{subfigure}{\smwdth\textwidth}
        \includegraphics[width=\linewidth]{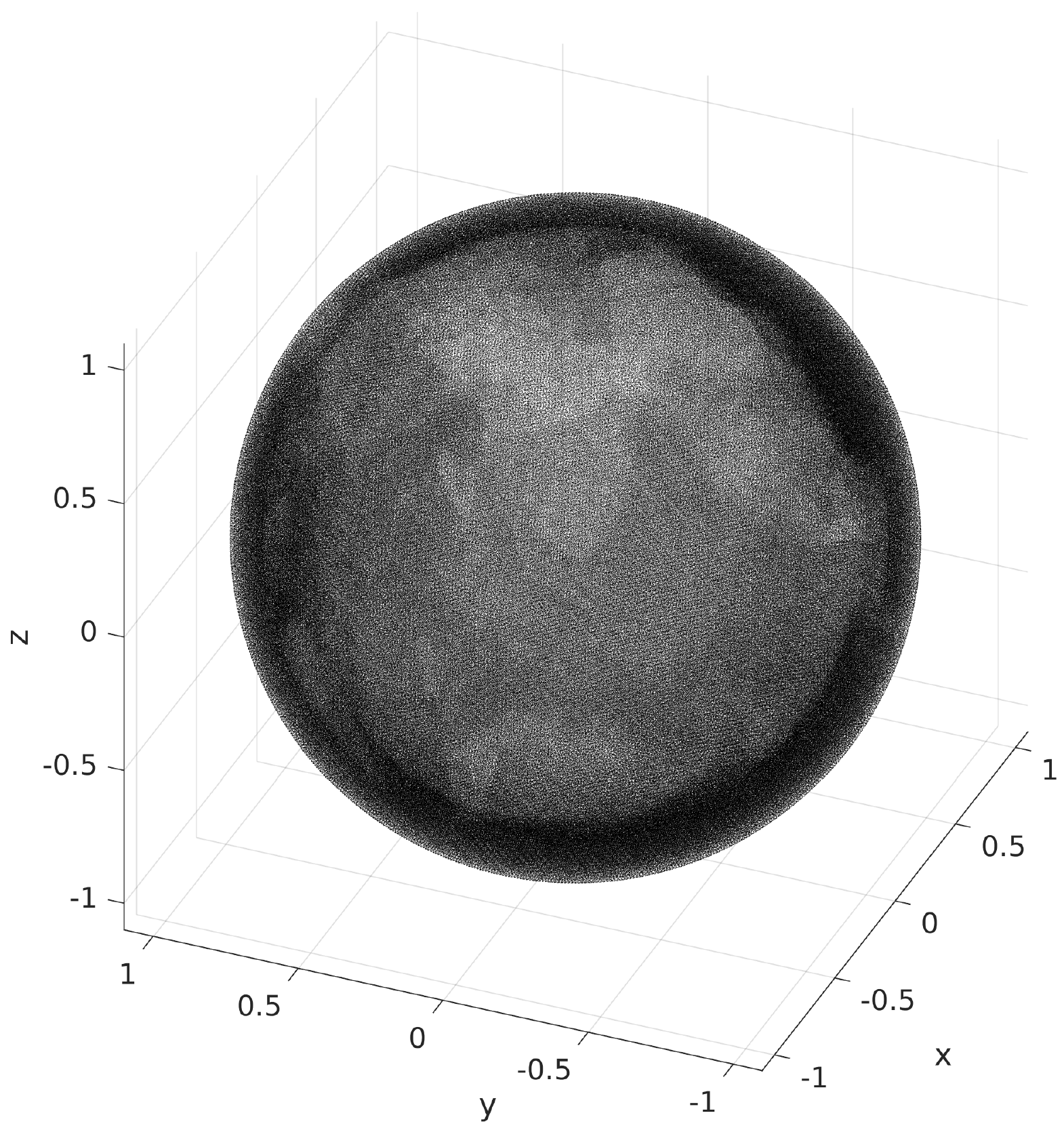} 
    \end{subfigure}
    ~\hspace*{.01\textwidth} 
    \begin{subfigure}{\smwdth\textwidth}
        \includegraphics[width=\linewidth]{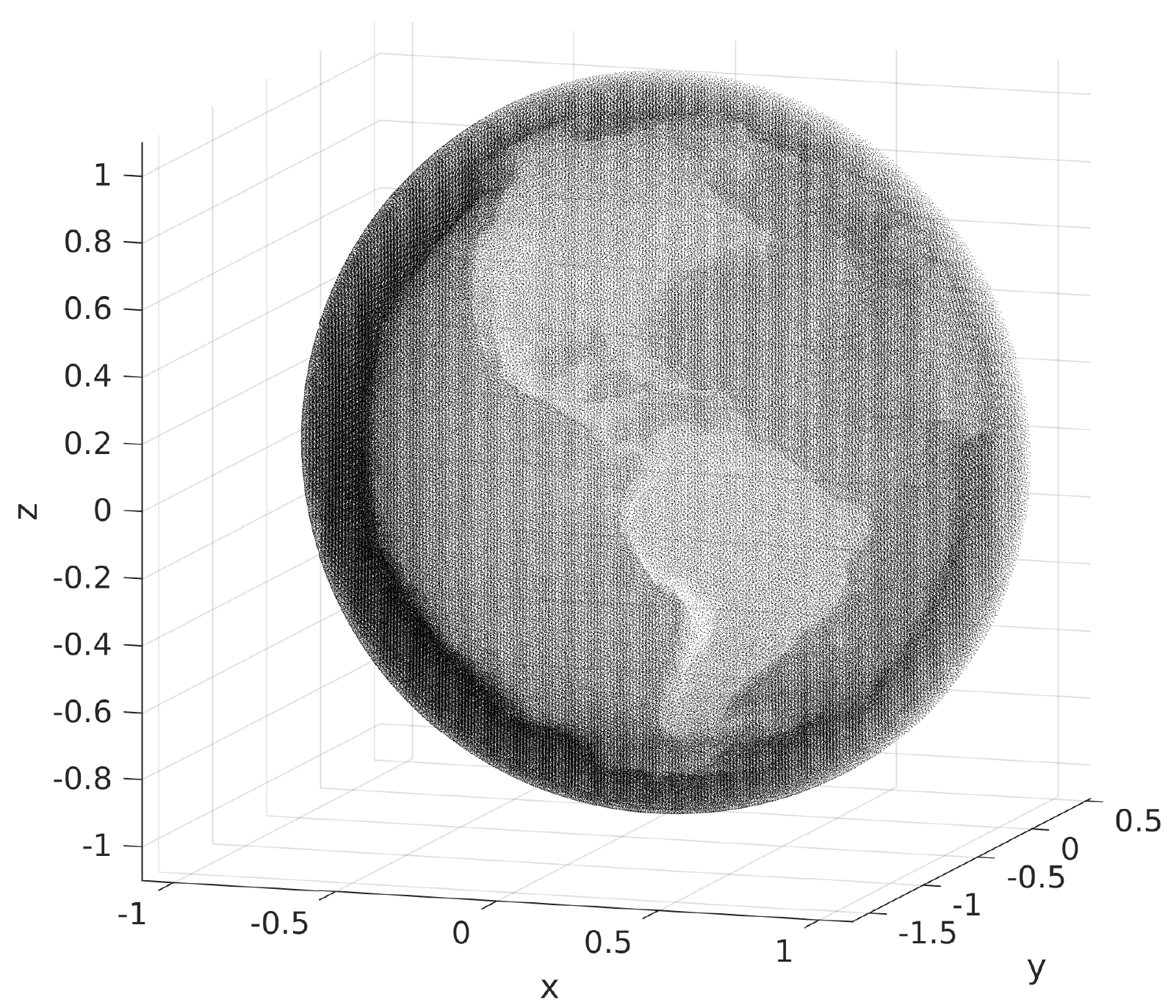} 
    \end{subfigure}
    \caption{Left: a general view of a uniform node distribution in an
    atmospheric-like shell. Right: a separate view of the Western hemisphere.} \label{fig:earth}
    \begin{subfigure}{\smwdth\textwidth}
        \includegraphics[width=\linewidth]{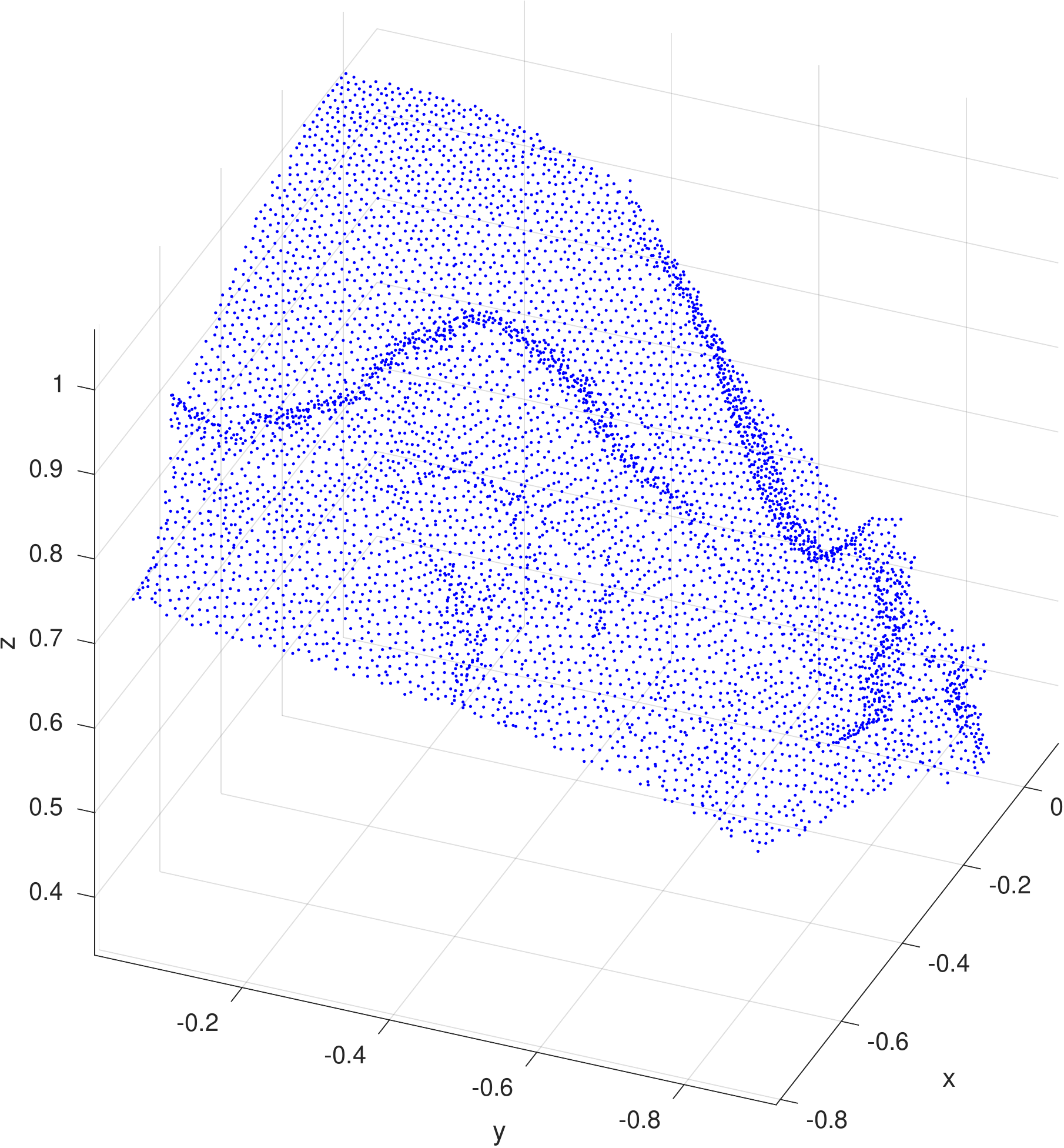} 
    \end{subfigure}
    ~\hspace*{.01\textwidth} 
    \begin{subfigure}{\smwdth\textwidth}
        \includegraphics[width=\linewidth]{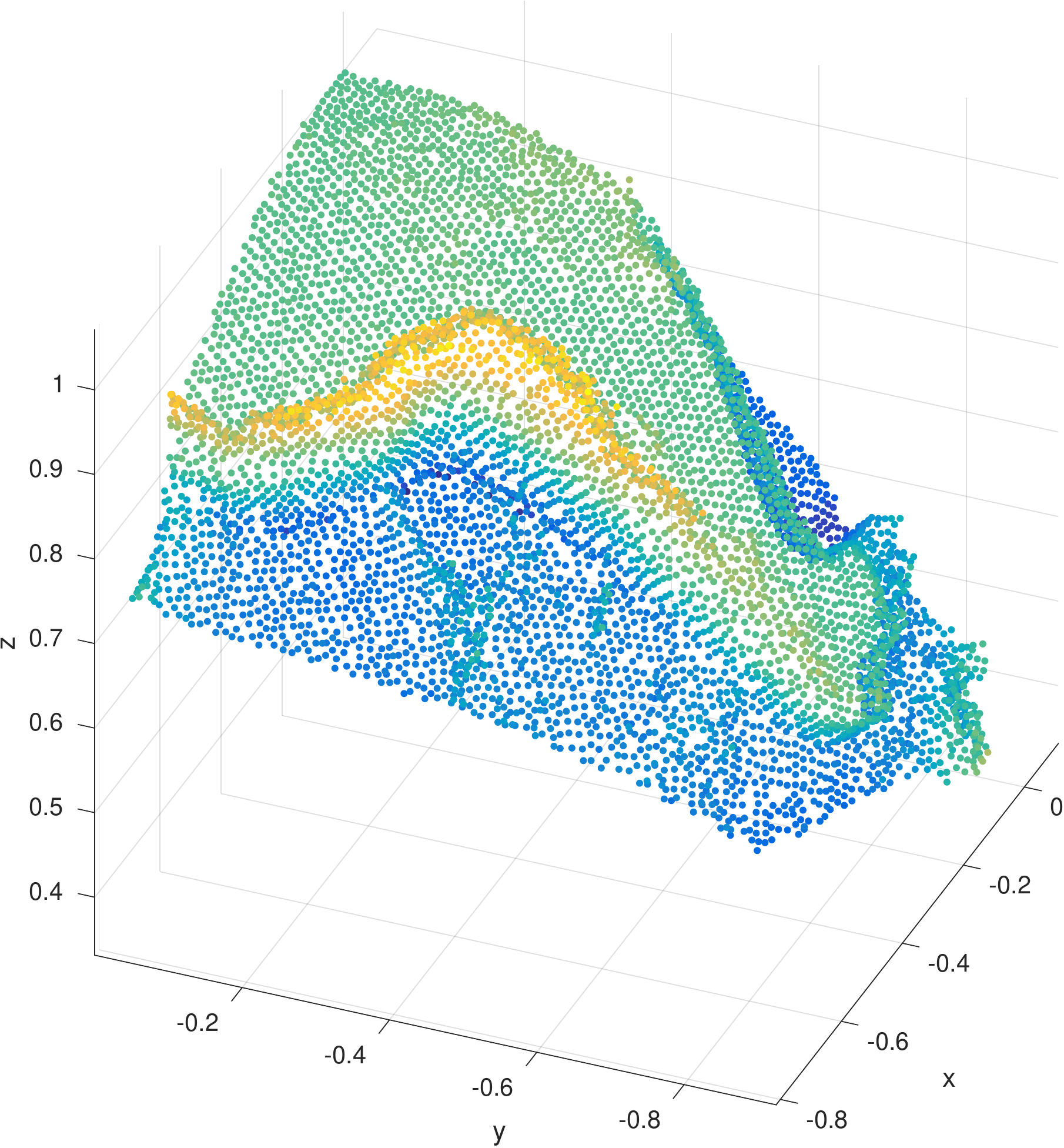} 
    \end{subfigure}
    \caption{Surface subset: a fragment of the Western coast of South America. The nodes on the right are color-coded using heights. } \label{fig:andes}
    \begin{subfigure}{\smwdth\textwidth}
        \includegraphics[width=\linewidth]{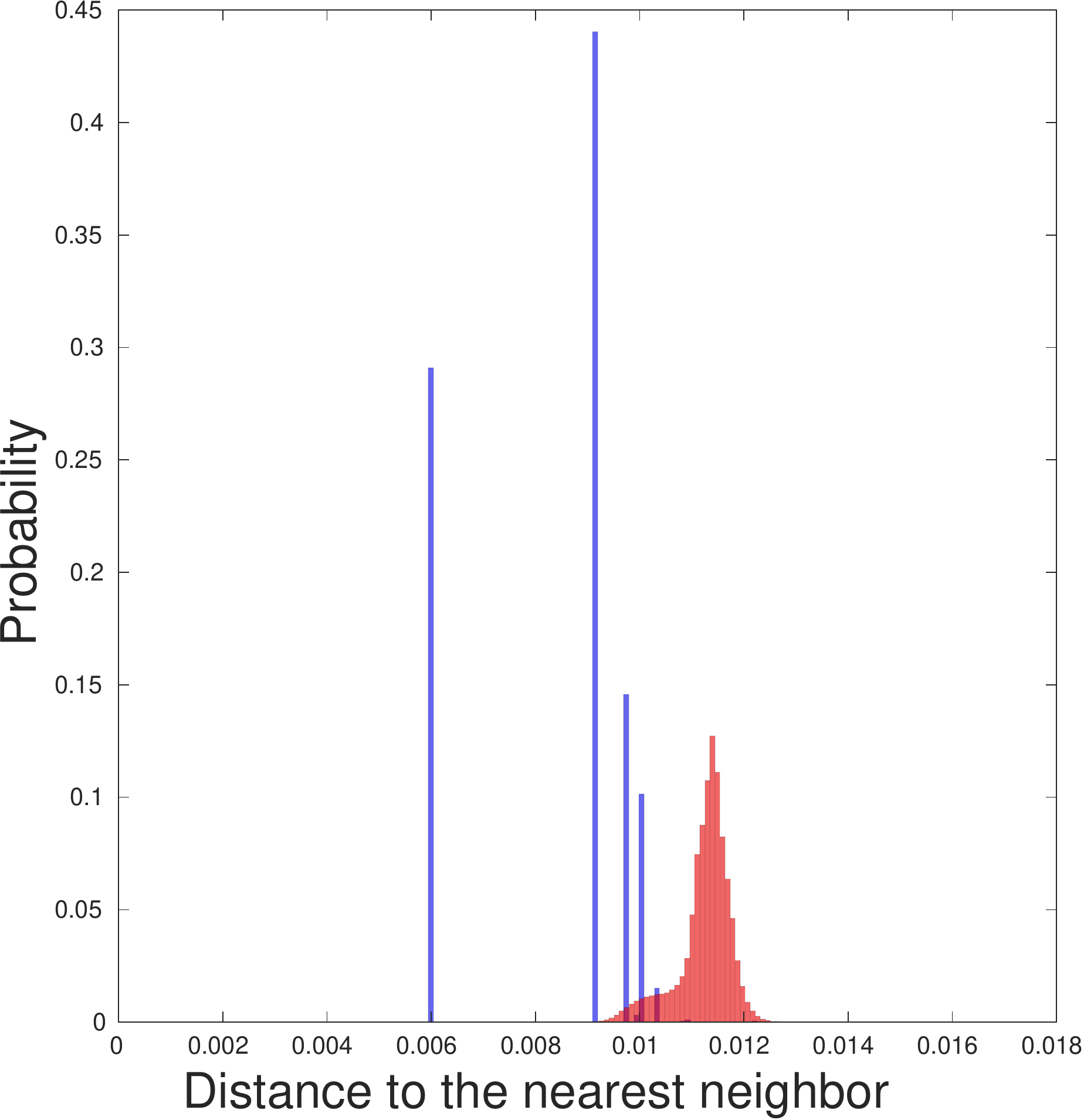} 
    \end{subfigure}
    ~\hspace*{.005\textwidth} 
    \begin{subfigure}{.40\textwidth}
        \includegraphics[width=\linewidth]{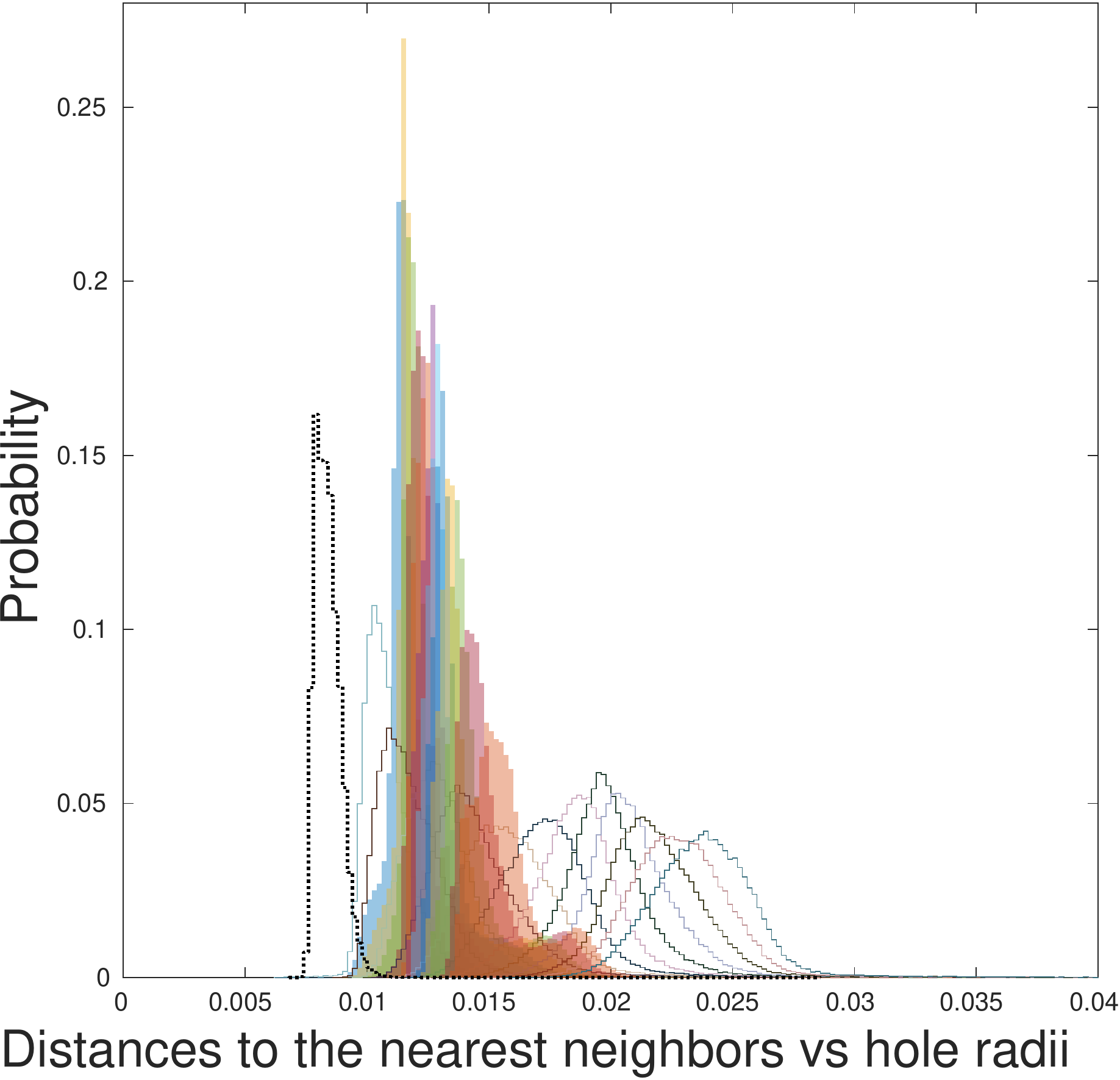}
    \end{subfigure}
    \caption{The effects of the repel procedure and hole radii.
        Left: probability distribution of the nearest-neighbor distances in the atmospheric node
        set, before (blue) and after (red) executing the repel subroutine. Right:
        distribution of distances to the 12 nearest neighbors for the whole
        configuration (color only), for the surface subset (contours), the hole radii
    (black dashed contour on the left).} \label{fig:hist_earth} 
\end{figure}
The ETOPO1 dataset stores relief as a 21\,600-by-10\,800 array of elevations above
the sea level;  equivalently, of radial coordinates that correspond to the
spherical angles defined by the array's indices. The data points are equispaced
on lines of constant azimuth/inclination with angular distance $ \AngleSpacing = {\pi}/{10\,800} $ between them. To determine whether a given node \label{p:spherical}$ \nd = (\VarRadius_{\nd},\VarAngle_{\nd},\VarPolarAngle_{\nd}) $ belongs to $ \Omega_\text{etopo} $, its radial coordinate $ \VarRadius_{\nd} $ was compared with a linear
interpolation of the values of radii of three ETOPO1 points with the nearest
spherical coordinates. For example, assume that such three points have the spherical
coordinates $ (\VarRadius_j,\VarAngle_j,\VarPolarAngle_j),\ j=1,2,3, $ where $ 0\leq\VarAngle\leq 2\pi $
and $ 0\leq\VarPolarAngle\leq\pi $ are the azimuth and polar angle, respectively, and
\[ \begin{array}{ccccc} 
        \VarAngle_1 =  l\AngleSpacing ,& 	\VarAngle_2 =  (l+1)\AngleSpacing,&
        \VarAngle_3 =  l\AngleSpacing, &0\leq l\leq 21\,599; \\ 
        \VarPolarAngle_1 =  m\AngleSpacing, & \VarPolarAngle_2 = m\AngleSpacing, & \VarPolarAngle_3 = (m+1)\AngleSpacing, 
& 0\leq m\leq 10\,799.\\ 
\end{array} \]
Without loss of generality, the inequalities 
\[ l\AngleSpacing \leq \VarAngle_{\nd} < (l+1)\AngleSpacing; \qquad  m\AngleSpacing \leq \VarPolarAngle_{\nd} < (m+1)\AngleSpacing\]
hold true.
The point inclusion function is defined in this case as  
\[
    \chi(\nd; {\Omega_\text{etopo}}) = 
    \begin{cases} 1, & \VarRadius_1 + \frac{\VarAngle_{\nd}-\VarAngle_1}{\AngleSpacing}(\VarRadius_2-\VarRadius_1) + \frac{\VarPolarAngle_{\nd}-\VarPolarAngle_1}{\AngleSpacing}(\VarRadius_3-\VarRadius_1) \,\bs{<}\, \VarRadius_{\nd} \,\bs{<}\, 1.1;\\
        0, & \text{otherwise},
    \end{cases}    
\]
with $ 1.1 $ being the radius of the outer sphere in the chosen scale.  In
effect, the algorithm for evaluating the $ \chi(\cdot\,;\Omega_\text{etopo}) $  described here coincides with the
star-shaped point location algorithm from \cite[Section 2.2]{Preparata1985}, applied to the interpolated Earth surface and the outer spherical boundary. 

Our node set consists of 1,356,566 nodes with the nearest-neighbor separation close to the constant $ \rho(\nd) = 0.01124 $, and our top priority was to ensure the low variance of the radial separation across the configuration, especially on the surface; the general view of the set is given in Figure~\ref{fig:earth}.  
We used the piecewise IL with golden-ratio derived parameters $\alpha_1 = \sqrt2,$  $ \alpha_2 = (\sqrt5-1)/\sqrt2$; regarding these $ \alpha_1,\, \alpha_2 $ see also the discussion in \hyperref[appendix]{Appendix}.
Several statistics of the resulting set are presented in the following table; here again we used the common notation $ \bar x $ for the averaged value of a quantity $ x $. Notation $ \Delta^k $ stands for the distance to the $ k $-th nearest neighbor.
\[
\begin{array}{lcc}
	\toprule 
    & \text{Whole node set} & \text{Surface nodes} \tabline \midrule 
	 \overline{\Delta^{12}(\nd) /\Delta^2(\nd)}    & 1.3674 & 2.0353  \tabline  
	 \overline{\Delta^4(\nd) /\Delta^1(\nd)}    & 1.0859  & 1.34019 \tabline  
	\text{99th percentile of } \{\Delta( \nd_i )\} & 0.012143 &0.014444  \tabline  
	 \overline{\Delta(\nd)}  &0.011243  & 0.010879 \tabline  
	\text{1st percentile of }\{\Delta( \nd_i )\} & 0.009652 & 0.009340 \tabline \bottomrule 
\end{array} 
\]

Figure~\ref{fig:andes} illustrates the distribution of nodes close to the
surface of $ \Omega_\text{etopo} $. No pullback function has been used, just the inclusion check performed as in \eqref{eq:node_update}.  Observe that the near-surface nodes display no artifacts, and the spacing does not significantly depend on the local surface shape.
The left subplot in the Figure~\ref{fig:hist_earth} illustrates the effect of
\ref{step:repel} on the distribution of distances to the nearest neighbor. In the right subplot, we have collected distances to the 12 nearest neighbors for the whole configuration, and separately for the surface subset.  
The histogram also
contains the distribution of \textit{hole radii}, that is, distances from the
Voronoi centers of the entire node configuration to their respective nearest
nodes. It is a well-known fact that the Voronoi centers are local maxima of the
distance from the node set \cite{Conway1999}, considered as a function on the
whole space $ \mathbb{R}^3 $.  Note that all the histograms on the right are
normalized by probability, not by the node count.

\begin{figure} \centering
    \includegraphics[width=\bgwdth\textwidth]{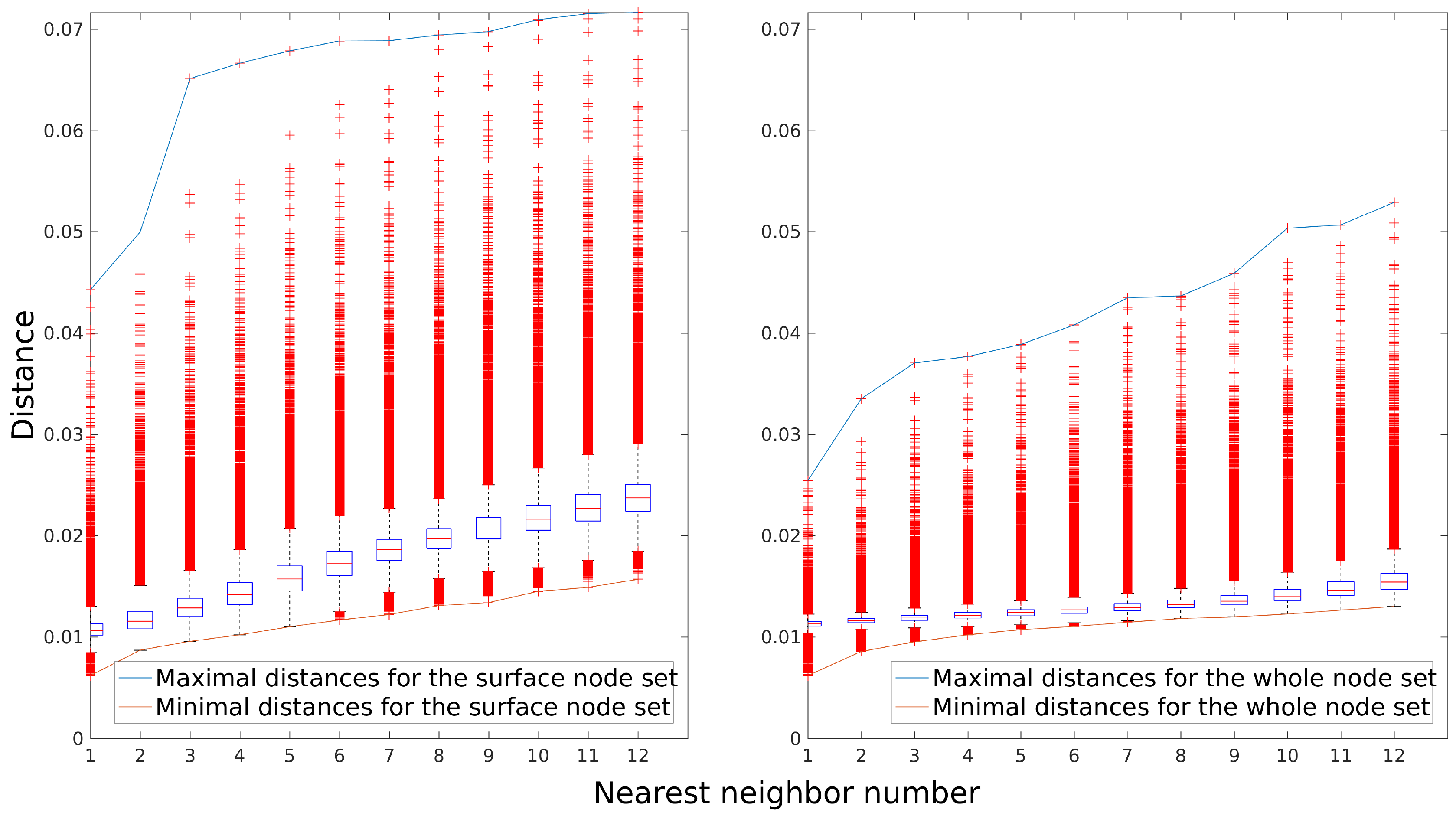}
    \caption{Distribution of distances to the 12 nearest neighbors for the
        atmospheric node configuration; medians and the 25th and 75th percentiles
        are shown. Left: the surface subset. Right: the whole set. Scales are the same
    in both subplots.} \label{fig:nns_earth} 
\end{figure}
The pair of plots in Figure~\ref{fig:nns_earth} shows in detail the
distribution of distances to the nearest neighbors in the sample node set. It
has been produced using the standard Matlab routine \textit{boxplot}. For each of
the blue boxes corresponding to a specific nearest neighbor, the central mark
is the median, the edges of the box denote the 25th and 75th percentiles. The
red crosses mark outliers.

\subsection{Point cloud}\label{subsec:cloud} To demonstrate a
nonuniform node distribution using our algorithm, we fix a collection of 100 points, $ \mathcal P_{100} $, inside the cube $ [-1,1]^3 $,  and consider the following radial density function: 
\[ \rho(\nd) = \left( \Delta(\nd;  \mathcal P_{100}  ) +\Delta^2(\nd;  \mathcal P_{100}  ) \right)/20, \]
where, as above, $ \Delta^k $  is the distance to the $ k $-th nearest neighbor. A possible interpretation of this density is a distribution that concentrates about a set of points $ \mathcal P_{100} $, which are of particular interest for a certain model. 

We proceed as in the \ref{step:repel} of the algorithm, not using the full gradient expression described in \eqref{eq:fullgradient}. In fact, it is instructive to note that computing the second term in \eqref{eq:fullgradient} would be quite cumbersome here in view of $ \rho $ being a piecewise-defined function. One could thus consider the density recovery for this distribution, Figure~\ref{fig:cloud_ratio}, as a validation of the gradient truncation approach in \ref{step:repel}; cf. Section~\ref{subsec:algo_discussed}.
The Q-MC voxels were drawn from the sequence $ \{ \lat_n\} $ with the same lattice parameters as in Section~\ref{subsec:atm_nodes}, $\alpha_1 = \sqrt2,$  $ \alpha_2 = (\sqrt5-1)/\sqrt2$.

\begin{figure} \centering 
    \begin{subfigure}{\smwdth\textwidth}
        \includegraphics[width=\linewidth]{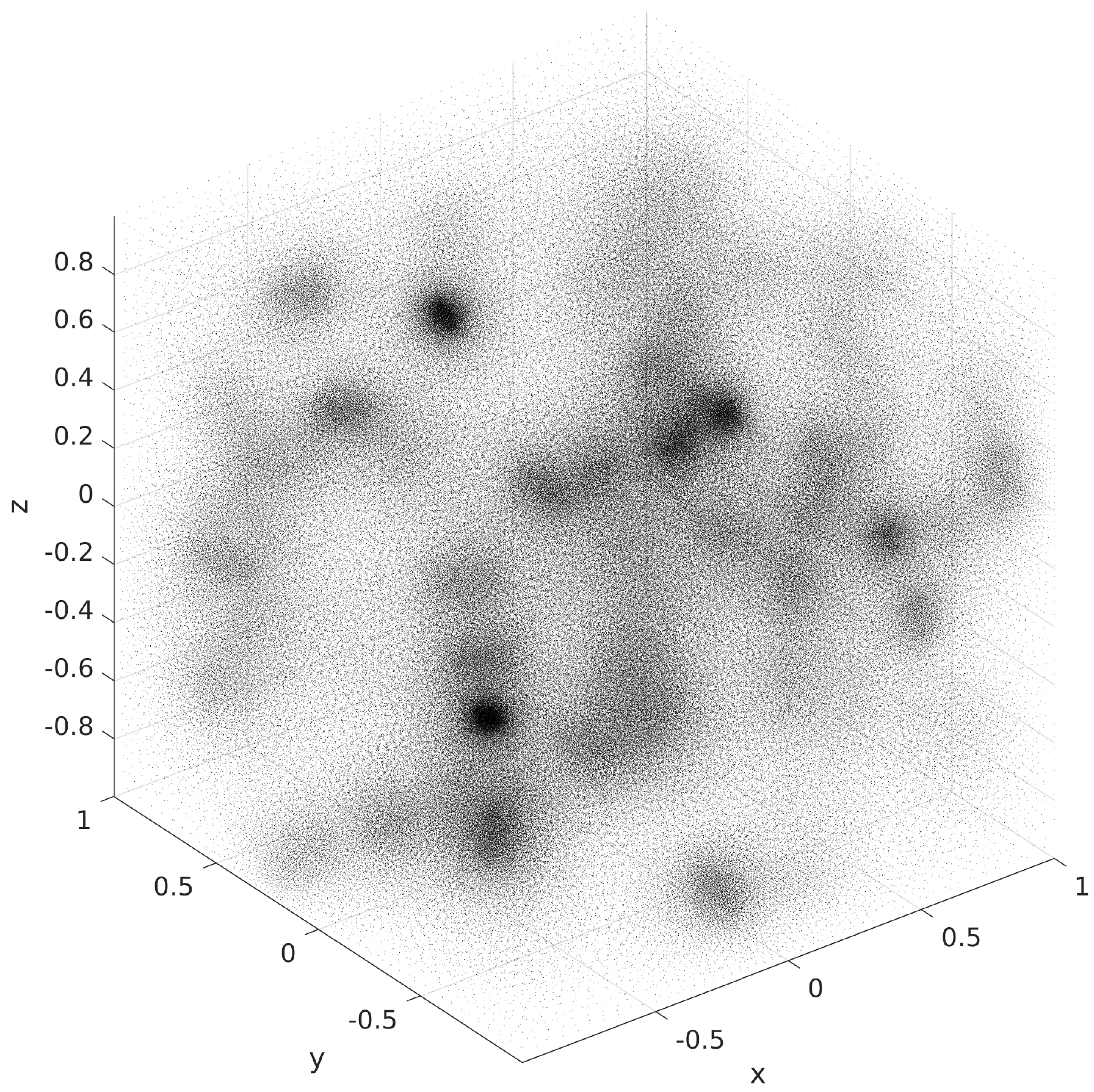} 
    \end{subfigure}
    ~\hspace*{.01\textwidth} 
    \begin{subfigure}{\smwdth\textwidth}
        \includegraphics[width=\linewidth]{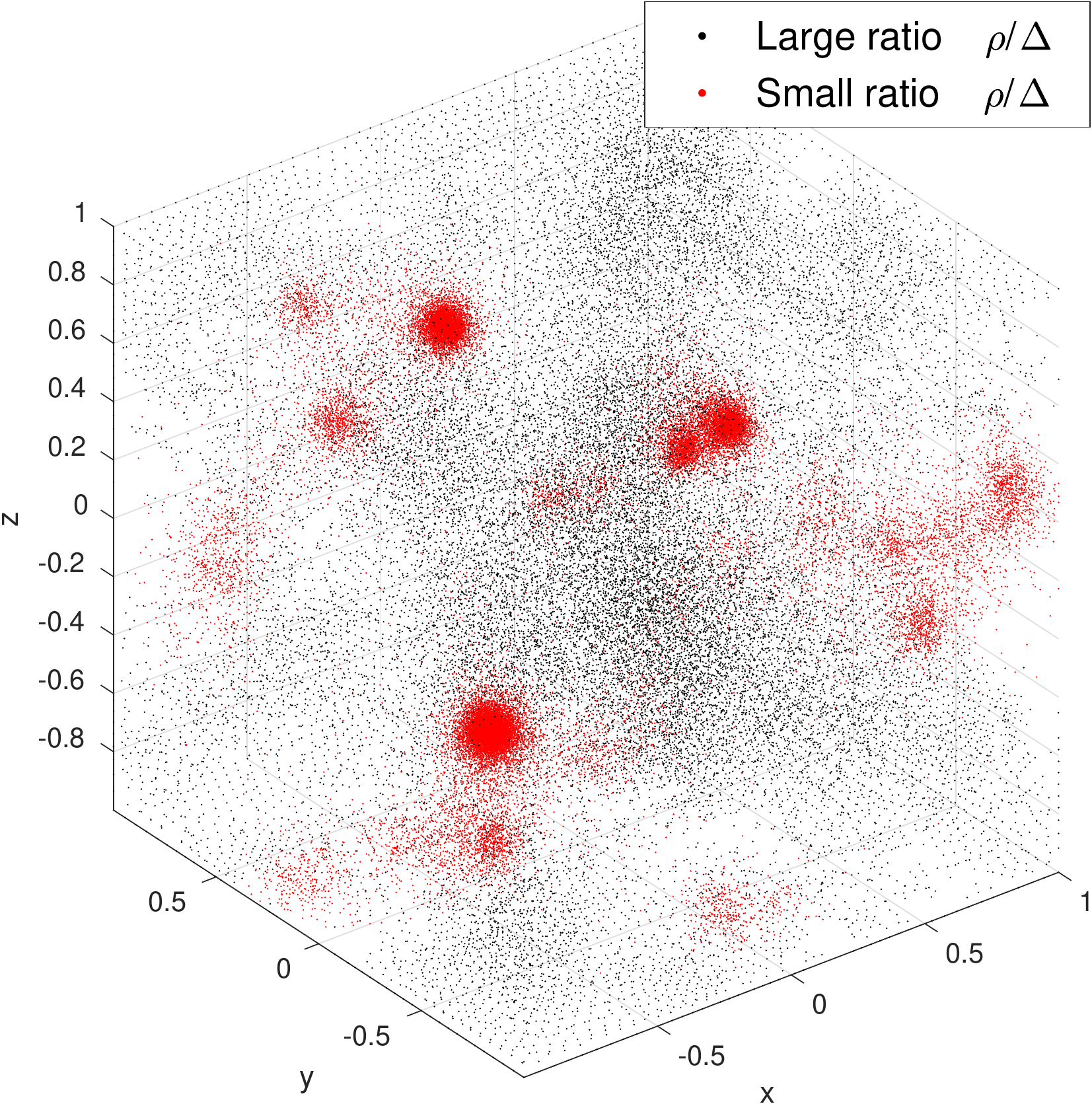}
    \end{subfigure} 
    \caption{Left: the node set from Section~\ref{subsec:cloud}.
    Right: node locations that contribute to the distribution of the ratio $ \rho(\nd)/\Delta(\nd) $ beyond the 5- and 95-percentiles.} \label{fig:cloud} 
\end{figure}

Figure~\ref{fig:cloud_ratio} contains the distribution of the ratio $
\rho(\nd)/\Delta(\nd)$. The minimal and maximal values of the ratio are about
0.8099 and 1.8231 respectively; its mean value is 0.9797, and the variance is
0.0019. The 5- and 95-percentiles are 0.9208 and 1.0441, respectively; the right plot in Figure~\ref{fig:cloud} highlights the outliers in the ratio distribution.
\begin{figure} \centering
    \includegraphics[width=\smwdth\textwidth]{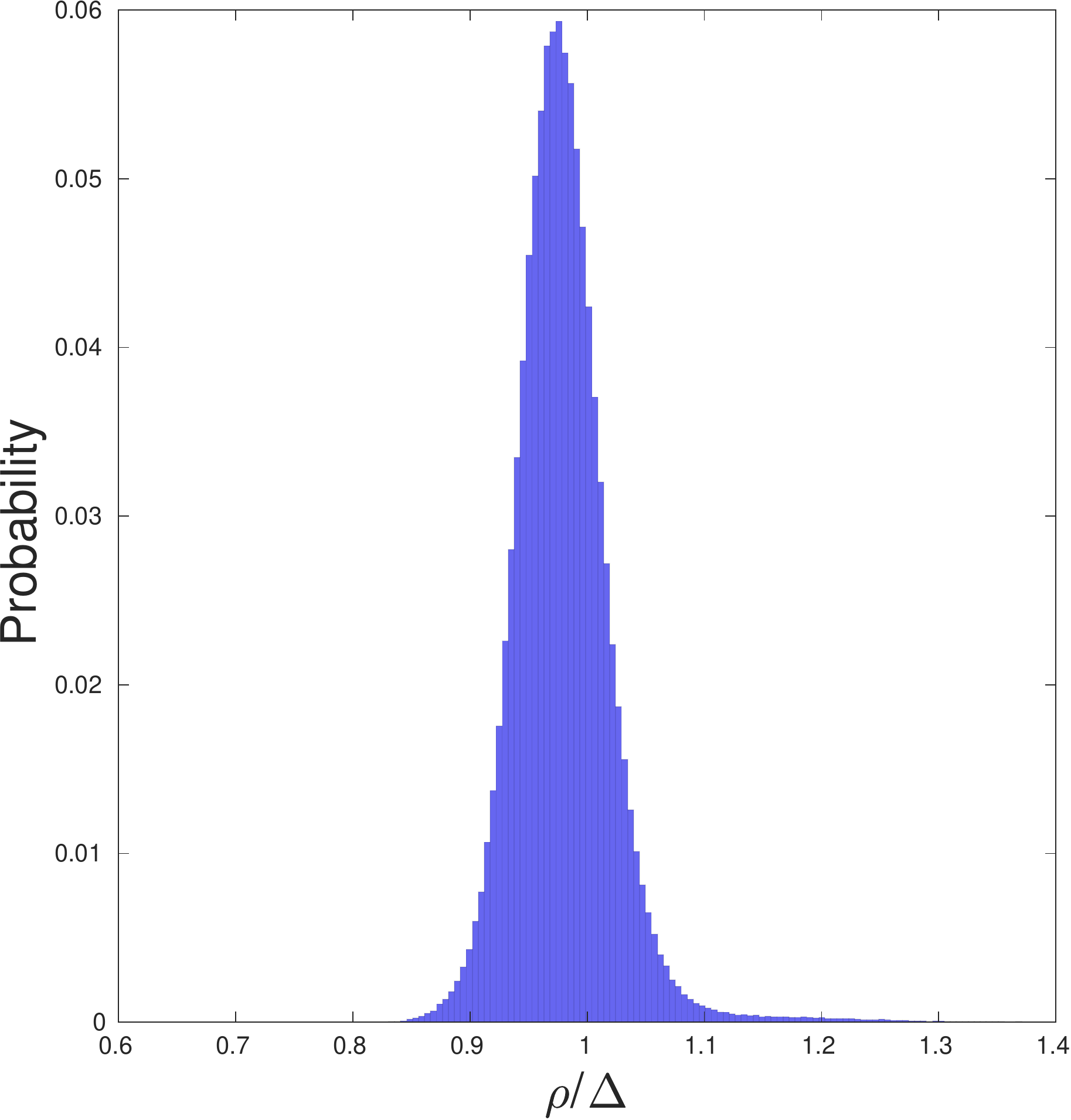}
\caption{Distribution of the ratios $ \rho(\nd)/\Delta(\nd) $ for the node set in Section~\ref{subsec:cloud}.} \label{fig:cloud_ratio} \end{figure}

\subsection{Spherical shell}\label{subsec:sph_shell} The motivation for this
example comes from atmospheric modeling. Representing the Earth surface
by a sphere, we consider first a thin 3-dimensional shell $ \Omega_\text{shell} $ of inner radius $R_\text{inner}$ and outer
radius $ R_\text{inner}+H_\text{atm}$ with constant target separation $h$ between points in the radial (vertical)
direction, and the tangential (horizontal) separation 
to be $\tau(r) = C\cdot r$ at radius $ r $, for some constant $C$. With typical choices of
parameters, $\tau$ will be much larger than $h$, reflecting the much higher
resolution needed in the vertical direction due to $ H_\text{atm} \ll R_\text{inner} $. We make a 
radial change of variables, which can written in spherical coordinates as $(r\,,\VarAngle\,,\VarPolarAngle) \mapsto (\hat{r}(r),\,\VarAngle,\,\VarPolarAngle)  $, so that any configuration in  $ \Omega_\text{shell} $ having the 2-directional resolutions $ \tau(r) $ and $ h $ will have isotropic resolution after the transformation. It is much easier to construct RBF bases in the isotropic case, hence our deliberation.

 Following this
change of variables, the radial/tangential node separations become, respectively, 
\begin{equation}\label{eq:2dir_separation}
\begin{aligned}
    \hat\nu(r) & = h\cdot \hat{r}'(r)\\
    \hat\tau(r) & = C\cdot \hat{r}(r).
\end{aligned}   
\end{equation}
Setting these two quantities to be equal, we obtain the ODE
\[\hat{r}'(r) = \frac{C}{h}\cdot\hat{r}(r)\]
with initial condition $\hat{r}(R_\text{inner}) = 1$, and its
solution becomes \[\hat{r}(r) = \exp\bigg(C\cdot\frac{r-R_\text{inner}}{h}\bigg).\]
From the second equation in \eqref{eq:2dir_separation} follows that our goal is to generate a node set in the $ (\hat{r},\,\VarAngle,\,\VarPolarAngle) $-space, whose separation is proportional to $\hat{r}$ and is equal in all directions: $ \rho(\nd) = C\cdot \|\nd\|$. The outer radius of the image of $ \Omega_\text{shell} $ in the $ (\hat{r},\,\VarAngle,\,\VarPolarAngle) $-space  is a function of $ R_\text{inner}  $ and $ H_\text{atm} $; our model implies $R_\text{inner} =$ 6,371,220, the mean radius of the Earth in meters, and $H_\text{atm} = $ 12,000, the thickness of the atmospheric layer we are interested in. The constant $ C $ is determined by the  desired tangential separation at the $ r=R_\text{inner} $ level.
\begin{figure} \centering 
    \begin{subfigure}{\smwdth\textwidth}
        \includegraphics[width=\linewidth]{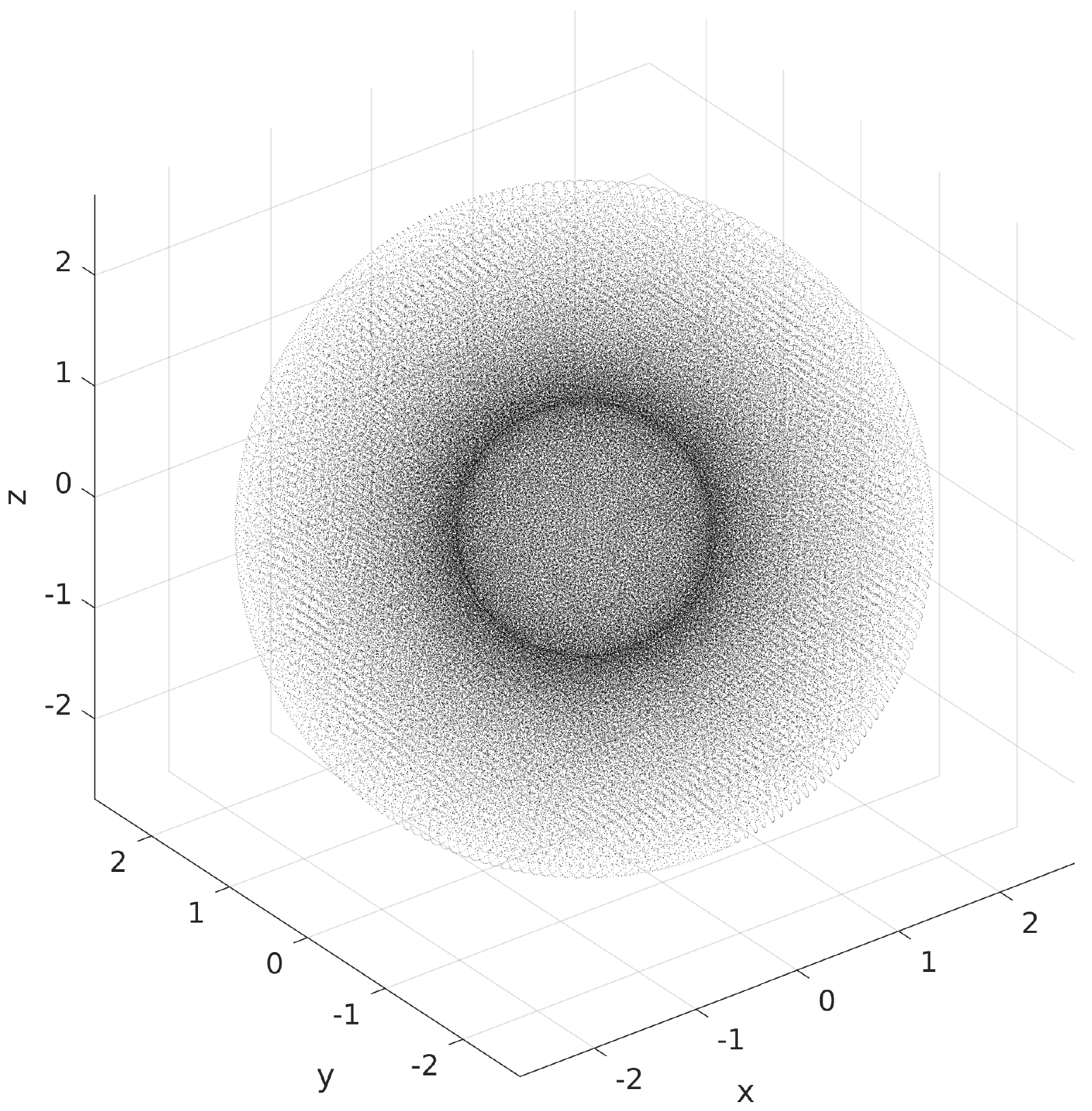} 
    \end{subfigure}
    ~\hspace*{.01\textwidth} 
    \begin{subfigure}{\smwdth\textwidth}
        \includegraphics[width=\linewidth]{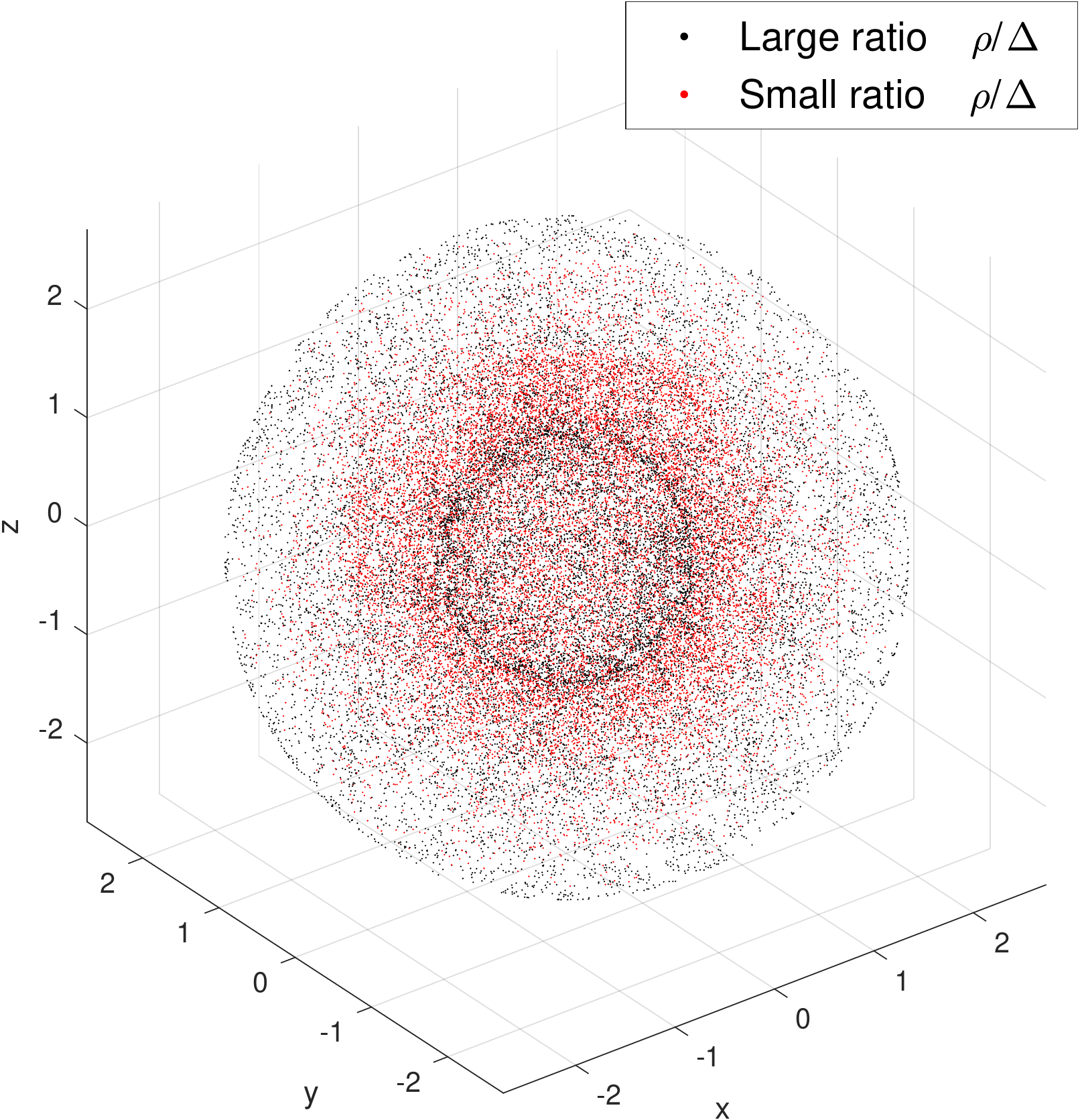}
    \end{subfigure} 
    \caption{Left: the node set from Section~\ref{subsec:sph_shell}.
    Right: node locations that contribute to the distribution of the ratio $ \rho(\nd)/\Delta(\nd) $ beyond the 5- and 95-percentiles.} \label{fig:sph_shell2} 
\end{figure}

\begin{figure} \centering 
    \begin{subfigure}{\smwdth\textwidth}
        \includegraphics[width=\linewidth]{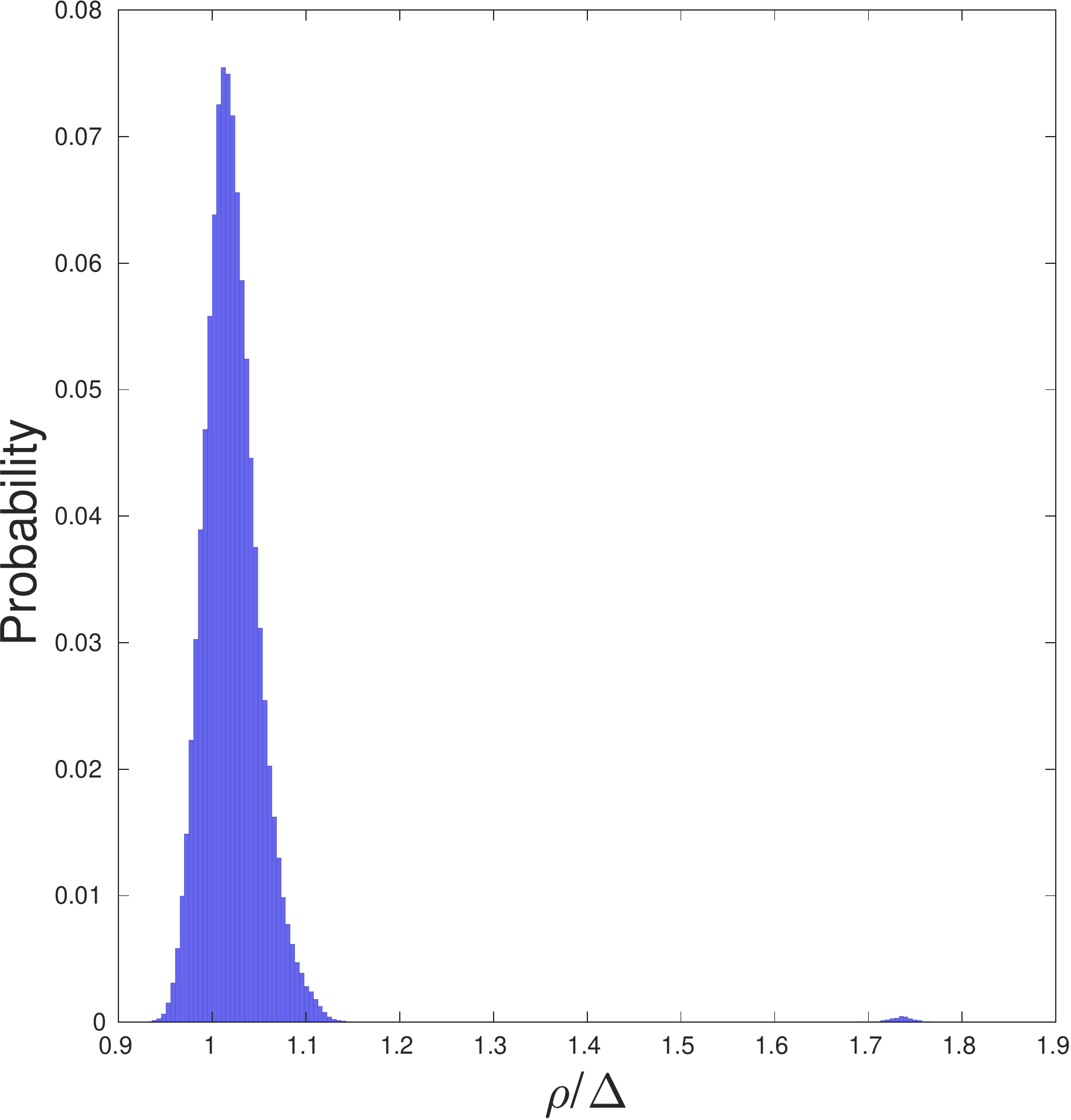}
    \end{subfigure} 
    ~\hspace*{.01\textwidth}
    \begin{subfigure}{\smwdth\textwidth} 
        \includegraphics[width=\linewidth]{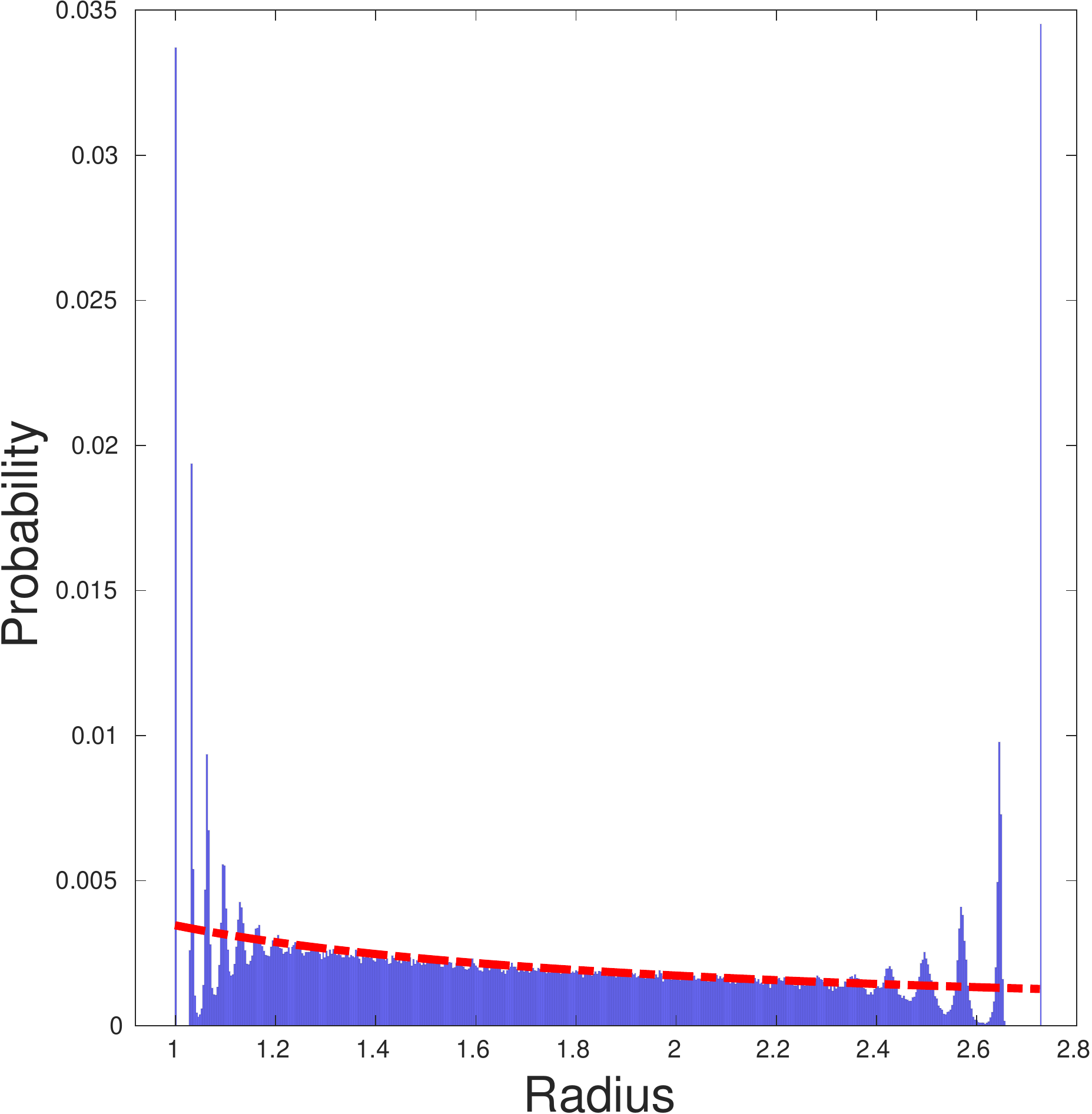} 
    \end{subfigure}
    \caption{Left: Distribution of the ratios $ \rho(\nd)/\Delta(\nd) $ for the node set in Section~\ref{subsec:sph_shell}. Right: Radial node distribution, actual (blue) and the theoretic (red) continuous component; i.e., without the $ \delta $-function spikes at the endpoints. }
    \label{fig:sph_shell} 
\end{figure}

Say, we intend to generate nodes corresponding to the 2 degree resolution on the spherical ``Earth surface" and $ h = 400$ meter vertical resolution. Due to the peculiarities of  atmospheric modeling, we would like to fix two much denser sets of nodes on the inner and outer boundary of $ \Omega_\text{shell} $; specifically, we are using 12,100 approximate Riesz energy minimizers on a sphere, appropriately rescaled. The interior nodes are generated using our algorithm, and then \ref{step:repel} is modified so as to leave the boundary subset intact. This, however, causes a difficulty: the separation distances between the interior and the surface nodes must remain large; on the other hand, our generic formulation of \ref{step:repel} does not account for the much higher concentration of nodes on the surface, which causes excessive repelling force, seen in the oscillations of the radial distribution, Figure~\ref{fig:sph_shell}. Mitigating this effect requires artificially weakening the repulsive force caused by the boundary nodes, a straightforward task using our codebase. Instead, we show in Figure~\ref{fig:sph_shell} the performance of the generic algorithm, to illustrate complications that may arise when applying it to specialized problems.

The set $ \Omega_\text{shell} $ can be challenging for the basic form of our algorithm, as described in Section~\ref{subsec:algo}: obtaining satisfying convergence requires using the full version of gradient descent \eqref{eq:fullgradient}. The reasons for it being more difficult to tackle than, say, $ \Omega_\text{etopo} $ in Section~\ref{subsec:atm_nodes}, are that due to convexity of the outer boundary, the weighted  $ s $-energy minimizers on it are denser than  on the sphere with radius $ \hat{r}(R_\text{inner}+H_\text{atm}) - 10^{-3} $, for example; see also discussion at the end of Section~\ref{subsec:riesz}. Getting rid of the artifacts at the endpoints of the radial distribution is done by using the full gradient, weakening the repulsion of the fixed boundary nodes, and \textit{not} striving for the full convergence of a minimization method applied to the Riesz energy.

In this example, we used the $ \{\rmin_n\} $ sequence to fill individual voxels. The left subplot in Figure~\ref{fig:sph_shell} contains the distribution of the ratio $ \rho(\nd)/\Delta(\nd)$. The minimal and maximal values of the ratio are about 0.9165 and 1.8989 respectively; its mean value is 1.0226, and the variance is 0.0024. The 5- and 95-percentiles are 0.9782 and 1.0717, respectively.

\subsection{Run times}\label{subsec:timings}
\begin{table}[h]
    \[
        \begin{array}{lccccc}
            \toprule 
            \text{Example}  &  K  &  T  & N &   \text{Q-MC distribution times, s}      & \text{Repel times, s}    \tabline \midrule 
            \text{Atmospheric nodes} &  33   & 29   &     1,356,566  &    5      & 89  \tabline  
            \text{Point cloud}  &  30   &  200 & 577,321 & 4  &  840 \tabline  
            \text{Spherical shell} &  30   &  200  &  358,915  &    1   &   144    \tabline \bottomrule
        \end{array} 
    \] 
    \caption{Timings of the examples in Sections~\ref{subsec:atm_nodes}--\ref{subsec:sph_shell}.}
    \label{tab:timing}
\end{table}
The execution times (in seconds) for the above examples are summarized in Table~\ref{tab:timing}, where, as before, $ K $ and $ T $ stand for the number of nearest neighbors and the number of iterations used in the repel procedure in \ref{step:repel}, respectively. The fifth column contains times required to fill the voxels selected at \ref{step:subcubes} with configurations from either $ \{\lat_n\} $ or $ \{\rmin_n\} $ and to remove any redundant nodes as in \ref{step:cleanup}. All the computations were performed on a dedicated machine with 40 GB RAM and an 8-core \textit{Intel Xeon} CPU.
Note that the basic Q-MC node sets for both sequences were precomputed, and the pre-computation times are not included in the table. Computation of configurations in $ \{\lat_n\} $ for $ 1\leq n \leq 200 $ took less than 1 second. An implementation of the $ \{\rmin_n\} $ sequence for $ 1 \leq n \leq 200 $ took 4311 seconds to generate; coordinates of the resulting minimizers as well as the corresponding average separation distances are distributed with the associated codebase \cite{VMFF17matlab}.  
\section{Final observations and comparisons}\label{sec:remarks}
\subsection{Comparisons}\label{subsec:comparisons} 
\begin{figure}
    \centering
    \begin{subfigure}{.38\textwidth}
        \includegraphics[width=\linewidth]{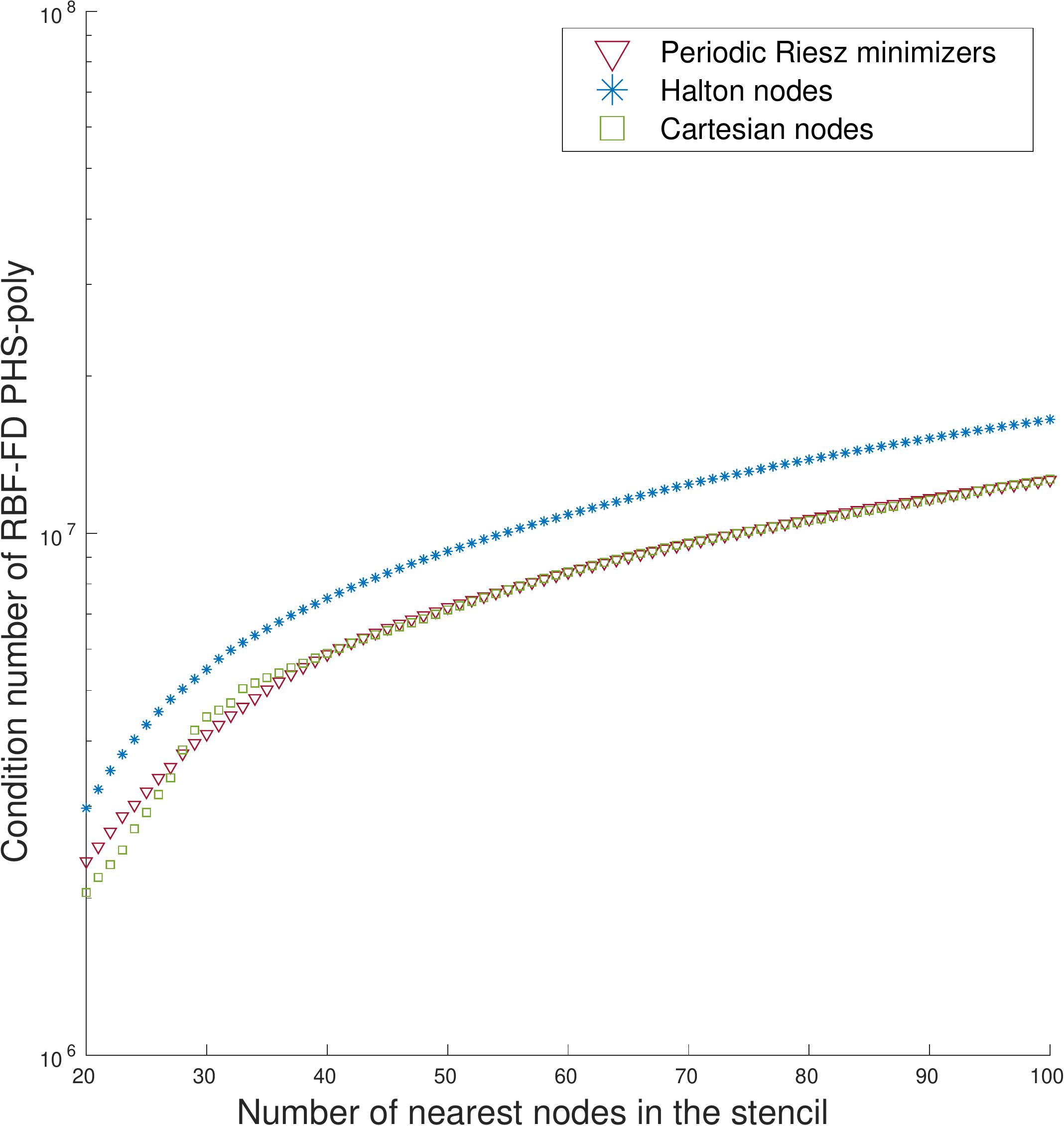} 
    \end{subfigure}
    ~\hspace*{.01\textwidth} 
    \begin{subfigure}{\smwdth\textwidth}
        \includegraphics[width=\linewidth]{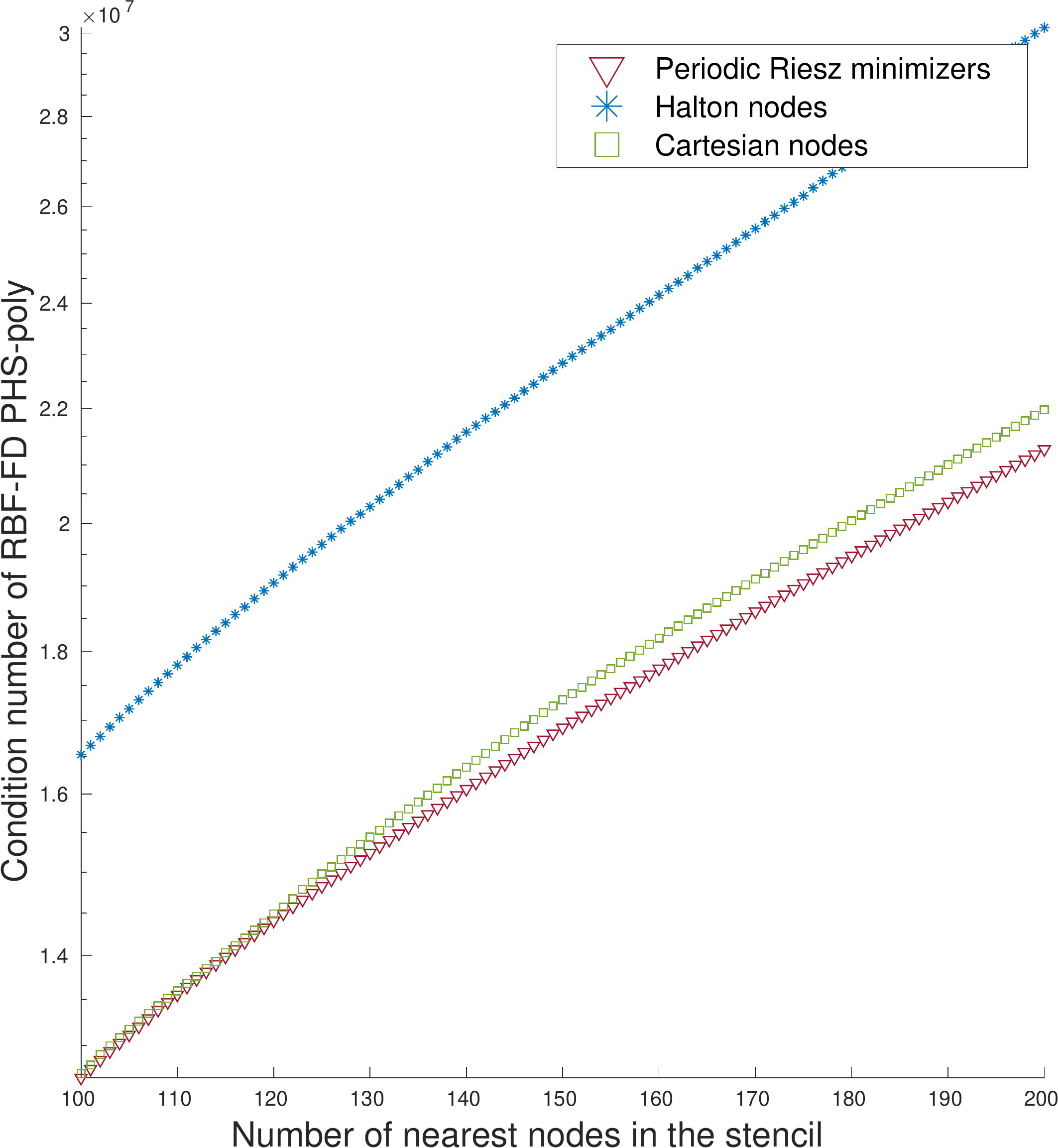} 
    \end{subfigure}
    \caption{ Average condition numbers of the joint RBF-FD PHS-based matrices for order 1- and 2-differential operators. {Left:} Between 20 and 100 nodes in the stencil.  {Right:} Between 100 and 200 nodes in the stencil.}
    \label{fig:condition_numbers} 
\end{figure}
Of the two Q-MC sequences we considered, the periodic Riesz minimizers appear more promising, being devoid of the lattice structure and having high space utilization. On the other hand, we have successfully applied ILs as an elementary uniform configuration.
One could use different sets of irrational parameters $ \alpha_1,\ldots,\alpha_{d-1} $ for different numbers of nodes in a voxel. Although this might be useful in mitigating the non-isotropic behavior of ILs, it makes hard to control node separation at voxel interfaces.  

Another quasi-uniform node set commonly used in Q-MC methods can be constructed from the Halton sequence \cite[Chapter 5.4]{Lemieux2009}, an example of a low-discrepancy sequence. To see how the Halton nodes compare to Riesz minimizers for RBF-FD methods, we have computed condition numbers of PHS-based RBF-FD matrix in the LHS of \eqref{eq:RBF-FD}, using operators $ \partial/\partial x_i,\, i=1,2,3 $, and $ \partial^2/\partial x_i \partial x_j,\, i\leq j = 1,2,3 $ as $ \mathfrak L $; the resulting LHS then constitutes the joint system for the weights corresponding to these nine differential operators. We used the RBF kernel $ \phi(r) = r^5 $, and the polynomials in the interpolation space were of degree at most $ 2 $, see Section~\ref{subsec:rbfs}. The computations were performed for the Riesz and Halton nodes, and the uniform Cartesian grid. The stencils consisted of $ K $ nearest neighbors of a random vector  with Gaussian distribution, centered around $ (0.5, 0.5, 0.5)^\text{tr} $; the evaluation point $ \bs x_0 $ was taken equal to the random vector itself. The nearest neighbors were drawn from 1000 nodes of the respective sequence, uniformly distributed over the unit cube. The Riesz nodes were produced by minimizing periodic energy \eqref{eq:s-riesz} with the distance \eqref{eq:pdist}.

Figure~\ref{fig:condition_numbers} contains a comparison of the condition numbers of RBF-FD matrices for the three sequences. The values shown are averages of the condition numbers for 500 random stencil centers $ \bs x_0 $; the averaging was introduced to eliminate the rather unpredictable dependence on $ \bs x_0 $, and to display the underlying trend. We omitted values of $ K $ below 20 from the plots, as all the three node sets resulted in relatively ill-conditioned systems; this was to be expected, as the recommended stencil size is roughly twice the number of linearly independent polynomials in the interpolation space (there are 10 monomials of degree at most 2 in $ \mathbb R^3 $) \cite{Flyer2016}.

\subsection{Range of applications}\label{subsec:applicability}
Our method has proven very efficient for slowly varying radial densities that are small (recall that small radial density means a large number of nodes per unit volume) compared to the entire node set scale, and is capable of handling very complex underlying sets.
The range of dimensions where the algorithm can be used efficiently is determined by the applicability of Q-MC initialization and the nearest neighbor searches: the repelling iterations for Riesz energy in \ref{step:repel} are largely (with a proper value of $ s $) dimension-agnostic. A shortcoming that is common to all quasi-Monte Carlo methods (but of little practical relevance) is a much worse performance (measured by $ L^2 $ discrepancy), compared to Monte-Carlo distribution, in dimensions starting at about $ 15 $ \cite{Caflisch1998}. Furthermore, using the uniform grid to detect the support $ \Omega $, as is done in \ref{step:subcubes}--\ref{step:densefill}, becomes unfeasible already for $ d = 10 $; instead, one needs an efficient way to determine which corners of the grid are in some sense close to $ \Omega $. This is certainly not a feature of our approach, but a manifestation of the curse of dimensionality: treatment of a complicated high-dimensional set is a computationally intensive task. Regarding finding nearest neighbors it should be noted that common implementations of k-d trees are efficient only up to about $ d=20$; additionally, the k-d tree approach is faster than the full brute force search only if $ N \gg 2^d $ \cite{Freidman1977}. On the other hand, as has already been noted, our repelling procedure does not require frequent updates of the search tree, as the local adjacency largely remains intact.

The suggested algorithm is very local, and it therefore must be straightforward to add multi-resolution and adaptive refinement, as is widely done for grids  \cite{clawpack, Debreu2008}, yet as of this writing, our proof-of-concept implementation does not include these features. Still, we would like to observe that refining the voxel structure is indeed easier than refining a mesh, since no geometry is taken into account. This partially addresses the previous remark on detection of $ \Omega $ in high dimensions.

The closest set of goals to what we have presented here, that we're aware of, is posed in the pioneering paper \cite{Shimada1995}; our method is crafted for full-dimensional domains, and apparently performs faster in this case. The bubble packing algorithm in \cite{Shimada1995} is conceptually similar to the greedy filling of centers in \ref{step:saturate}, while physical relaxation is an alternative to the energy minimization we employ; of course, the idea of relaxation can also be found in a number of related references, and is a well-known approach in this context, see for example \cite{Persson2004a}.
Our method requires computing the gradient of the desired radial density in the cases when the outer boundary of the underlying set is uniformly convex, and/or when the radial density changes quickly. Alternatively, fine partition of the set is necessary. Either solution, however, may be computationally expensive. 

    \section{Appendix: separation properties of sequences \texorpdfstring{ $ \{\lat_n\}$ and $ \{\rmin_n\} $}{\{Ln\} and  \{Mn\}} }\label{appendix} 
This Appendix deals with the results of our numerical experiments, set in the $3$-dimensional space. 
\begin{figure}
    \centering
    \begin{subfigure}{\fwdth\textwidth}
        \includegraphics[width=\linewidth]{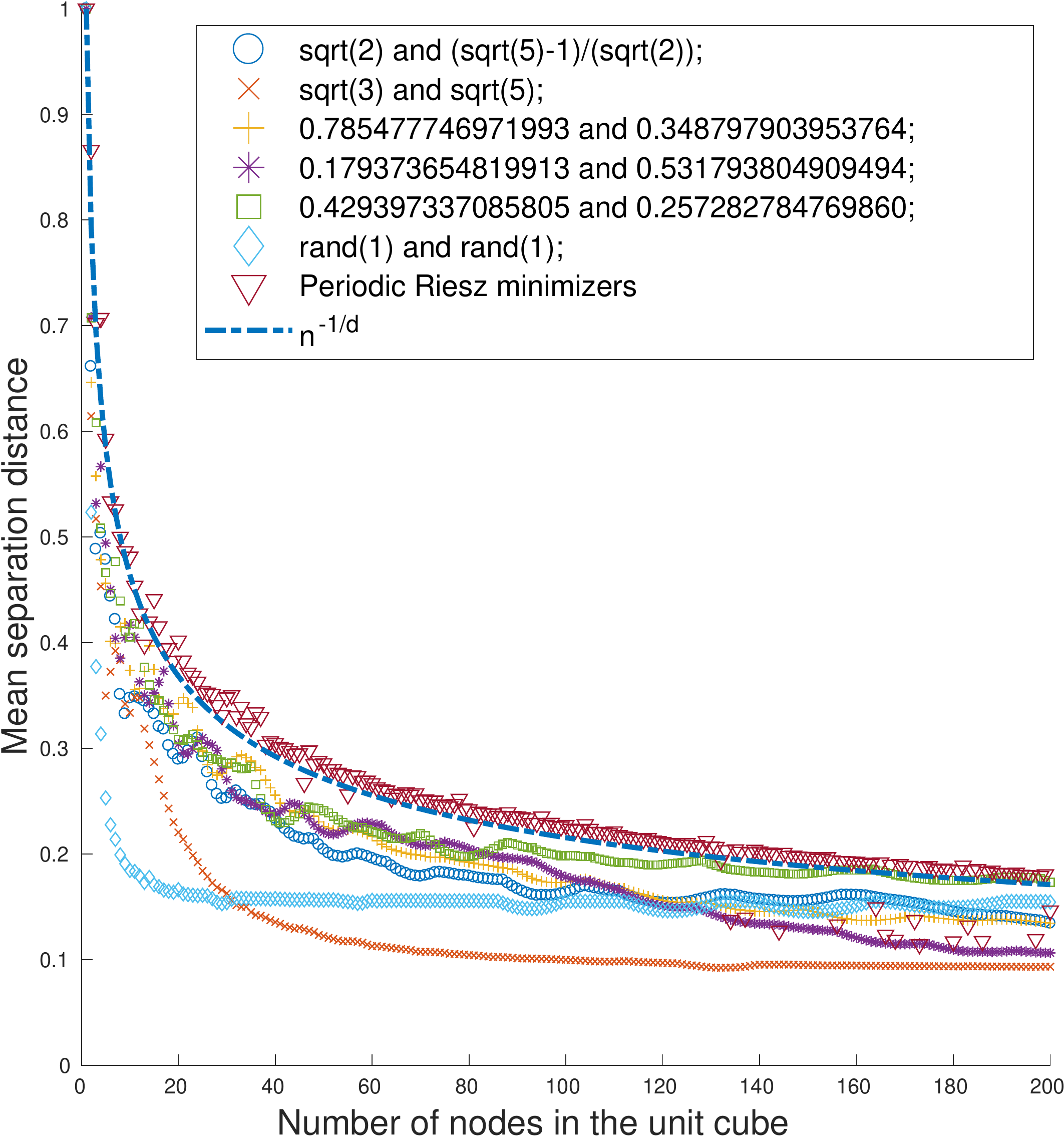} 
    \end{subfigure}
    ~\hspace*{.01\textwidth} 
    \begin{subfigure}{\fwdth\textwidth}
        \includegraphics[width=\linewidth]{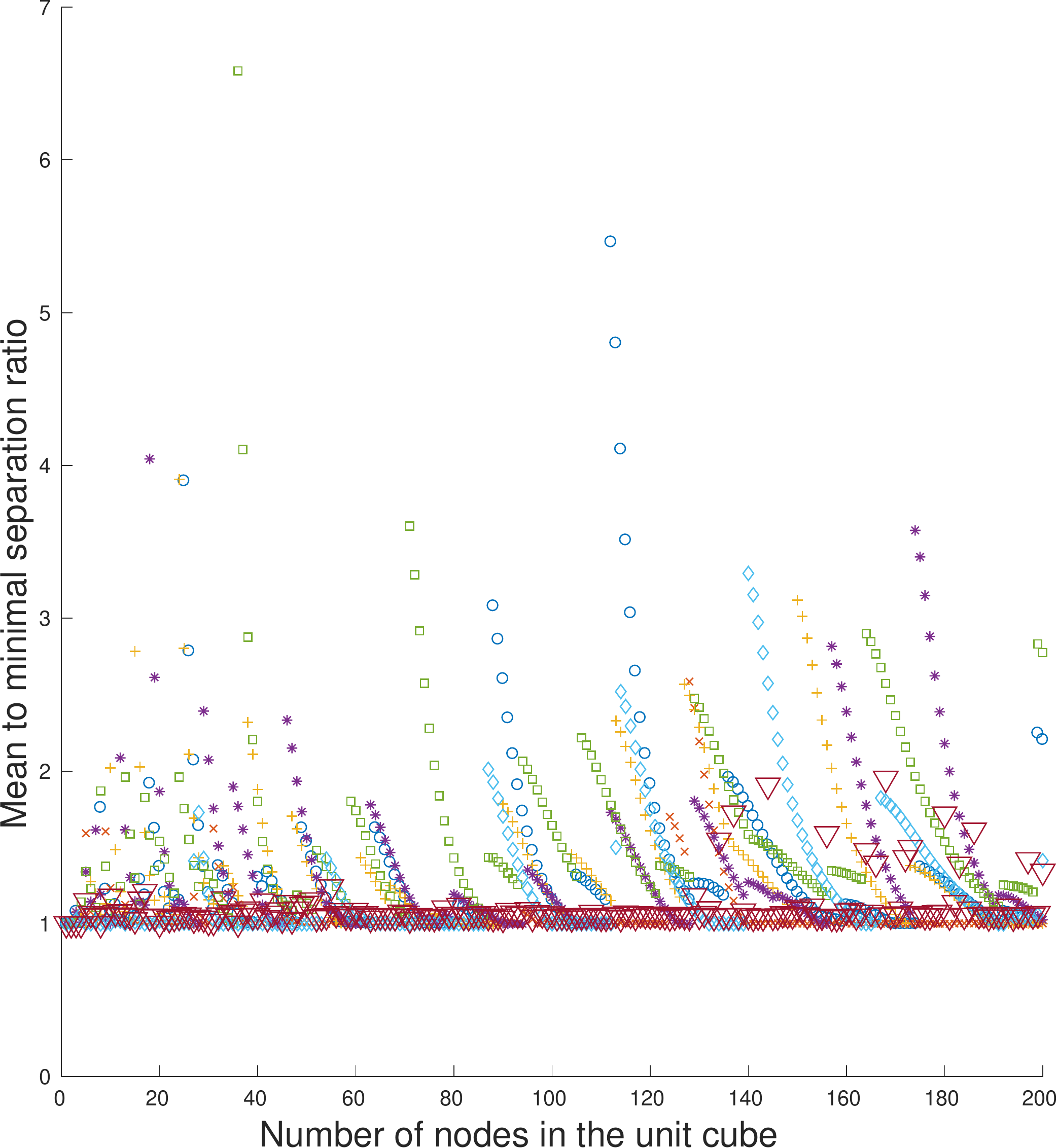} 
    \end{subfigure}
    \caption{{Left:} dependence of the mean separation distances on the number of nodes in the unit cube for different values of parameters $ \alpha_1,\ \alpha_2 $; the $ n^{-1/d} $ decay rate shown as a dash-dot line.  {Right:} ratios of the mean separation distances to the minimal ones for the same configurations. }
    \label{fig:lattices} 
\end{figure}
The function $ \lambda(r) $ used in \ref{step:densefill} is the
number of nodes in the unit cube $ [0,1]^3 $, placed according to \eqref{eq:irr_lat}, or
obtained by minimizing the Riesz $ s $-energy \eqref{eq:s-riesz} with periodic metric, such that the
mean separation distance of these nodes is the closest to $ r $. To compute $ \lambda(r) $ for the periodization of $ \{\lat_n\} $, we tabulate mean separations $ \bar\Delta_n $ in a sample configuration
comprising $ \lat_n $ and its  $ 26 = 3^3-1 $ copies, obtained translating $ \lat_n $ by the vectors  $ \bs{\{}(i,j,k)^\text{tr} \bs{:} i,j,k \in \{0,\pm1\} \text{ and } |i|+|j|+|k|>0 \bs{\}}  $. The 
tabulated dependence of separation on $ n $ is then inverted and interpolated using a piecewise cubic Hermite interpolating polynomial. The reason to consider
separation distance between configurations in $ 3^d $ cubes in dimension $ d $ (and not a single cube with a single instance of $ \lat_n $) is to account for the boundary effects. Likewise, to compute $ \lambda(r) $ for the Riesz minimizers, the mean separation of $ \rmin_n $ is tabulated for $ 1\leq n \leq n_\text{max} $, then the inverse dependence is interpolated. No copies of $ \rmin_n $ are considered alongside the original configuration, since periodicity condition is already included in the metric \eqref{eq:pdist}.

In general, putting too many nodes in individual voxels is justified only if
the radial density function $ \rho $ varies slowly. For our applications, $
n_\text{max} \leq 100 $ was sufficient. 
\begin{figure}
    \centering
    \begin{subfigure}{\smwdth\textwidth}
        \includegraphics[width=\linewidth]{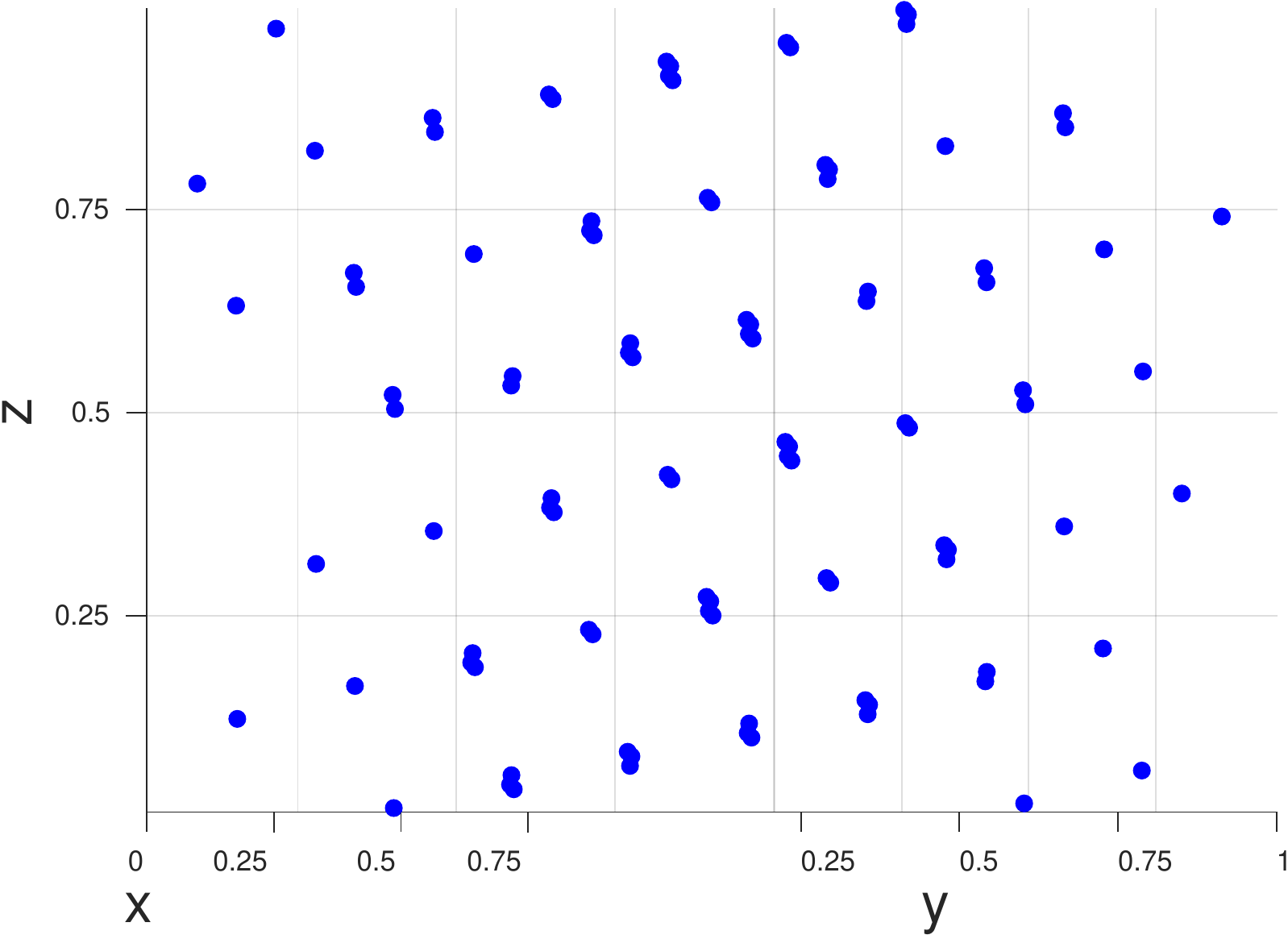} 
    \end{subfigure}
    ~\hspace*{.01\textwidth} 
    \begin{subfigure}{\smwdth\textwidth}
        \includegraphics[width=\linewidth]{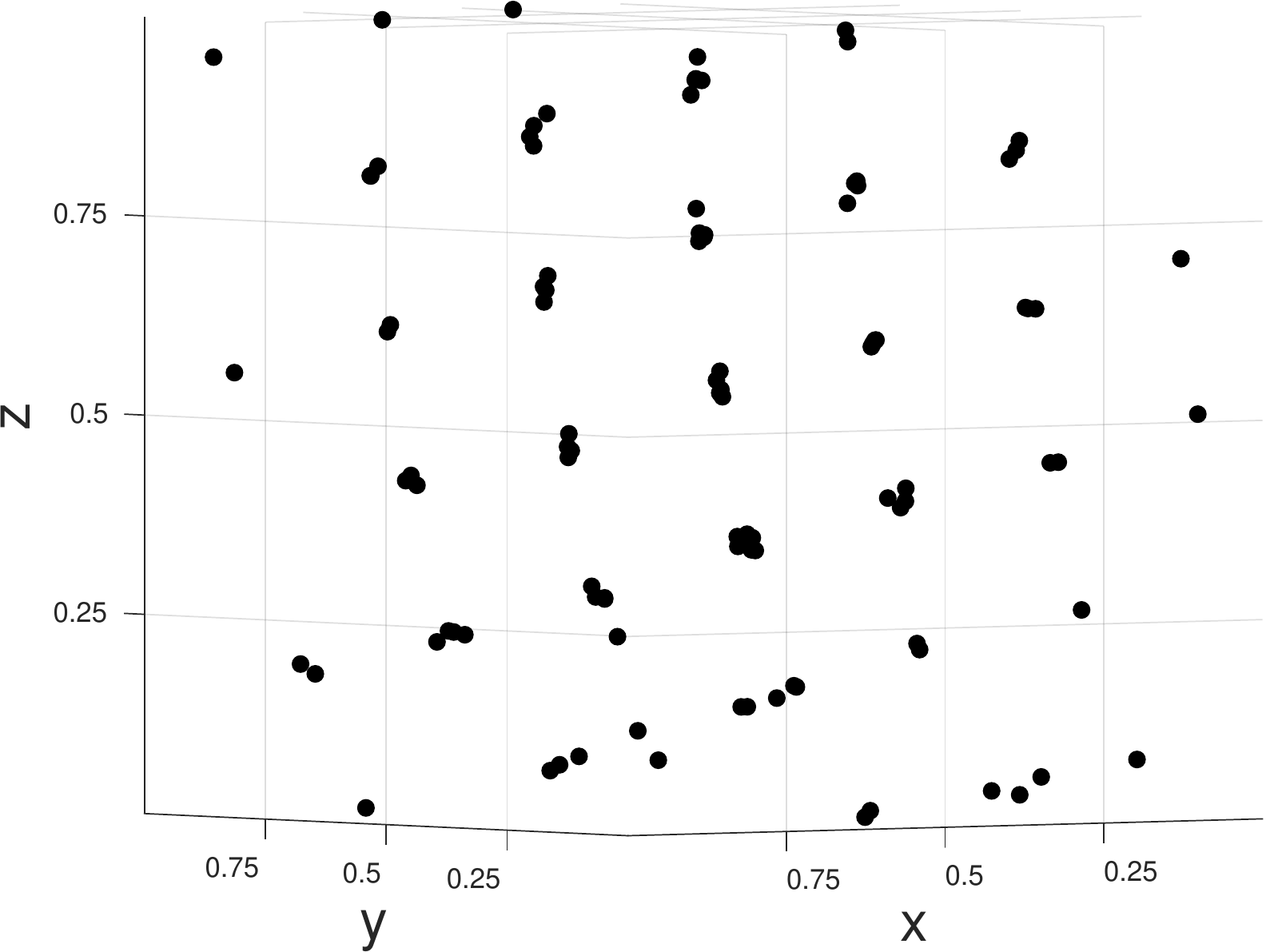} 
    \end{subfigure}
    \caption{{Left:} A cross-section of the IL $ \lat_{100} $ with parameters $ 0.179373654819913 $ and $  0.531793804909494 $. {Right:} A (different) cross-section of $ \rmin_{100} $. }
    \label{fig:cross-sections} 
\end{figure}
The left plot in Figure~\ref{fig:lattices} illustrates the delicate
dependence of the separation distances of ILs on the lattice parameters. While any set of irrational quantities $\alpha_1, \ldots,
\alpha_{d-1}$ in \eqref{eq:irr_lat} that are linearly independent over rationals
will give a uniformly distributed IL as $n$ grows, certain values  may
perform better than the others. In particular, adjustments can be made to
improve the distribution for small values of $n$. For example, it is known from \cite{Bilyk2013} that  a
2-dimensional IL generated by the golden ratio has optimal  $ L^2 $
discrepancy. 
Numerical experiments have shown that its 3-dimensional analog with parameters $\alpha_1 = \sqrt2,\,	 \alpha_2 = (\sqrt5-1)/\sqrt2$ does perform well for large numbers of nodes; yet by carrying out a Monte Carlo search for the parameters
maximizing separation distance in \eqref{eq:irr_lat}, we found several (necessarily rational) pairs that performed at least just as well for up to $ n= 200$, see Figure~\ref{fig:lattices}.

Curiously enough, a pair of random numbers drawn uniformly from $ [0,1] $ (shown in the legend as rand(1)), consistently performed better than the pair $ \sqrt3 $ and $ \sqrt5 $, starting at $ n \approx 40 $. We were able to reproduce this behavior in a number of runs; in fact, we haven't seen a random pair that wouldn't always outperform $ \sqrt3 $ and $ \sqrt5 $ after a fairly small $ n $.

The second graph in Figure~\ref{fig:lattices} shows the ratios of the mean  to minimal separation distances $ \bar\Delta_n/\Delta_n $ for the same range of $ n $. In both subfigures, Riesz periodic minimizers clearly stand out, by having the largest mean separation (left), and by smallest ratios (right). This means, the nearest neighbor distances $ \Delta(\nd) $ vary little from node to node in the $ \{\rmin_n\} $ sequence. We conclude this section by presenting in Figure~\ref{fig:cross-sections} a pair of cross-sections of the IL $ \lat_{100} $ and the configuration $ \rmin_{100} $ that look remarkably similar. In fact, we found the vague resemblance between the low-energy periodic configurations and lattice structures, similar to ILs, quite interesting, given the connection between packing and Riesz energy minimization \cite{Hardin2005}, and that the highest packing density in the 3-dimensional space is achieved, in particular, by the hcp lattice \cite{Conway1999}.

\section*{Acknowledgements} We thank Edward Saff and Douglas Hardin for their interest and useful comments. O.V. and T.M. were supported, in part, by the NSF
grant DMS-1516400; they would also like to express their gratitude to the University of Colorado, Boulder, where a part of this work was completed. The National Center for Atmospheric
Research is sponsored by the National Science Foundation.

\bibliography{nodes} \bibliographystyle{acm}

\begin{thebibliography}{10}

\bibitem{Alishahi2015a}
{\sc Alishahi, K., and Zamani, M.}
\newblock {The spherical ensemble and uniform distribution of points on the
  sphere}.
\newblock {\em Electron. J. Probab. 20\/} (2015).

\bibitem{ETOPO1}
{\sc Amante, C., and Eakins, B.~W.}
\newblock {ETOPO1 1 Arc-Minute Global Relief Model: Procedures, Data Sources
  and Analysis}.
\newblock \url{https://doi.org/10.7289/V5C8276M}, 2009.

\bibitem{Bayona2017}
{\sc Bayona, V., Flyer, N., Fornberg, B., and Barnett, G.~A.}
\newblock {On the role of polynomials in RBF-FD approximations: II. Numerical
  solution of elliptic PDEs}.
\newblock {\em J. Comput. Phys. 332\/} (2017), 257--273.

\bibitem{Beltran2016a}
{\sc Beltr{\'{a}}n, C., Marzo, J., and Ortega-Cerd{\`{a}}, J.}
\newblock {Energy and discrepancy of rotationally invariant determinantal point
  processes in high dimensional spheres}.
\newblock {\em J. Complex. 37\/} (2016), 76--109.

\bibitem{Bilyk2013}
{\sc Bilyk, D.}
\newblock {The $ L^2 $ discrepancy of irrational lattices}.
\newblock In {\em {Monte Carlo and Quasi-Monte Carlo Methods 2012}}. Springer,
  2013, pp.~289--296.

\bibitem{Bilyk2012}
{\sc Bilyk, D., Temlyakov, V.~N., and Yu, R.}
\newblock {The $ L^2 $ Discrepancy of Two-Dimensional Lattices}.
\newblock In {\em Springer Proc. Math. Stat.}, vol.~25. 2012, pp.~63--77.

\bibitem{Bollig2012}
{\sc Bollig, E.~F., Flyer, N., and Erlebacher, G.}
\newblock {Solution to PDEs using radial basis function finite-differences
  (RBF-FD) on multiple GPUs}.
\newblock {\em J. Comput. Phys. 231}, 21 (2012), 7133--7151.

\bibitem{Borodachov2008}
{\sc Borodachov, S.~V., Hardin, D.~P., and Saff, E.~B.}
\newblock Asymptotics for discrete weighted minimal {R}iesz energy problems on
  rectifiable sets.
\newblock {\em Trans. Amer. Math. Soc. 360}, 3 (2008), 1559--1580.

\bibitem{Borodachov2014}
{\sc Borodachov, S.~V., Hardin, D.~P., and Saff, E.~B.}
\newblock {Low Complexity Methods For Discretizing Manifolds Via Riesz Energy
  Minimization}.
\newblock {\em Found. Comput. Math. 14}, 6 (2014), 1173--1208.

\bibitem{Brauchart2012b}
{\sc Brauchart, J.~S., Hardin, D.~P., and Saff, E.~B.}
\newblock {The next-order term for optimal Riesz and logarithmic energy
  asymptotics on the sphere}.
\newblock {\em Contemp. Math 578\/} (2012), 1--31.

\bibitem{Brauchart2016}
{\sc Brauchart, J.~S., Reznikov, A.~B., Saff, E.~B., Sloan, I.~H., Wang, Y.~G.,
  and Womersley, R.~S.}
\newblock Random point sets on the sphere{\textemdash}hole radii, covering, and
  separation.
\newblock {\em Experimental Mathematics\/} (2016), 1--20.

\bibitem{BroomheadD.S.andLowe1988}
{\sc {Broomhead, D. S. and Lowe}, D.}
\newblock {Multivariable Functional Interpolation and Adaptive Networks}.
\newblock {\em Complex Syst. 2\/} (1988), 321-- 355.

\bibitem{Buhmann2003}
{\sc Buhmann, M.~D.}
\newblock {\em {Radial Basis Functions}}.
\newblock Cambridge University Press, Cambridge, 2003.

\bibitem{Buhmann2010}
{\sc Buhmann, M.~D., Dinew, S., and Larsson, E.}
\newblock {A note on radial basis function interpolant limits}.
\newblock {\em IMA J. Numer. Anal. 30}, 2 (2010), 543--554.

\bibitem{Caflisch1998}
{\sc Caflisch, R.~E.}
\newblock {Monte Carlo and quasi-Monte Carlo methods}.
\newblock {\em Acta Numerica 7\/} (1998), 1.

\bibitem{Chang2011}
{\sc Chang, C., and Lin, C.}
\newblock {LIBSVM}.
\newblock {\em ACM Trans. Intell. Syst. Technol. 2}, 3 (2011), 1--27.

\bibitem{MR2474372}
{\sc Cheney, W., and Light, W.}
\newblock {\em A course in approximation theory}, vol.~101 of {\em Graduate
  Studies in Mathematics}.
\newblock American Mathematical Society, Providence, RI, 2009.
\newblock Reprint of the 2000 original.

\bibitem{clawpack}
{\sc {Clawpack Development Team}}.
\newblock Clawpack software, 2017.
\newblock Version 5.4.0.

\bibitem{Conway1999}
{\sc Conway, J.~H., and Sloane, N. J.~A.}
\newblock {\em Sphere packings, lattices and groups}, third~ed., vol.~290 of
  {\em Grundlehren der Mathematischen Wissenschaften [Fundamental Principles of
  Mathematical Sciences]}.
\newblock Springer-Verlag, New York, 1999.
\newblock With additional contributions by E.~Bannai, R.~E.~Borcherds,
  J.~Leech, S.~P.~Norton, A.~M.~Odlyzko, R.~A. Parker, L.~Queen and
  B.~B.~Venkov.

\bibitem{Debreu2008}
{\sc Debreu, L., Vouland, C., and Blayo, E.}
\newblock {AGRIF}: Adaptive grid refinement in fortran.
\newblock {\em Computers {\&} Geosciences 34}, 1 (2008), 8--13.

\bibitem{MR2357267}
{\sc Fasshauer, G.~E.}
\newblock {\em Meshfree approximation methods with {MATLAB}}, vol.~6 of {\em
  Interdisciplinary Mathematical Sciences}.
\newblock World Scientific Publishing Co. Pte. Ltd., Hackensack, NJ, 2007.
\newblock With 1 CD-ROM (Windows, Macintosh and UNIX).

\bibitem{Flyer2016a}
{\sc Flyer, N., Barnett, G.~A., and Wicker, L.~J.}
\newblock {Enhancing finite differences with radial basis functions:
  Experiments on the Navier-Stokes equations}.
\newblock {\em J. Comput. Phys. 316\/} (jul 2016), 39--62.

\bibitem{Flyer2016}
{\sc Flyer, N., Fornberg, B., Bayona, V., and Barnett, G.~A.}
\newblock {On the role of polynomials in RB-FD approximations: I. Interpolation
  and accuracy}.
\newblock {\em J. Comput. Phys. 321\/} (2016), 21--38.

\bibitem{Flyer2012}
{\sc Flyer, N., Lehto, E., Blaise, S., Wright, G.~B., and St-Cyr, A.}
\newblock {A guide to RBF-generated finite differences for nonlinear transport:
  Shallow water simulations on a sphere}.
\newblock {\em J. Comput. Phys. 231}, 11 (2012), 4078--4095.

\bibitem{Flyer2013}
{\sc Flyer, N., Wright, G.~B., and Fornberg, B.}
\newblock {Radial Basis Function-Generated Finite Differences: A Mesh-Free
  Method for Computational Geosciences}.
\newblock {\em Handb. Geomathematics\/} (2013), 1--30.

\bibitem{Fornberg1996}
{\sc Fornberg, B.}
\newblock {\em A practical guide to pseudospectral methods}.
\newblock Cambridge University Press, 1996.

\bibitem{Fornberg1998b}
{\sc Fornberg, B.}
\newblock {Calculation of Weights in Finite Difference Formulas}.
\newblock {\em SIAM Rev. 40}, 3 (1998), 685--691.

\bibitem{Fornberg2015b}
{\sc Fornberg, B., and Flyer, N.}
\newblock {\em {A Primer on Radial Basis Functions with Applications to the
  Geosciences}}.
\newblock Society for Industrial and Applied Mathematics, Philadelphia, PA,
  2015.

\bibitem{Fornberg2015a}
{\sc Fornberg, B., and Flyer, N.}
\newblock {Fast generation of 2-D node distributions for mesh-free PDE
  discretizations}.
\newblock {\em Comput. Math. with Appl. 69}, 7 (2015), 531--544.

\bibitem{Fornberg2015}
{\sc Fornberg, B., and Flyer, N.}
\newblock Solving {PDEs} with radial basis functions.
\newblock {\em Acta Numerica 24\/} (2015), 215--258.

\bibitem{FornbergLarssonFlyer}
{\sc Fornberg, B., Larsson, E., and Flyer, N.}
\newblock {Stable Computations with Gaussian Radial Basis Functions}.
\newblock {\em SIAM J. Sci. Comput. 33}, 2 (jan 2011), 869--892.

\bibitem{Fornberg2013}
{\sc Fornberg, B., Lehto, E., and Powell, C.}
\newblock {Stable calculation of Gaussian-based RBF-FD stencils}.
\newblock {\em Comput. Math. with Appl. 65}, 4 (feb 2013), 627--637.

\bibitem{Fornberg2004b}
{\sc Fornberg, B., Wright, G., and Larsson, E.}
\newblock {Some Observations Regarding Interpolants in the Limit of Flat Radial
  Basis Functions}.
\newblock {\em Comput. Math. with Appl. 47}, 1 (2004), 37--55.

\bibitem{Fornberg2007}
{\sc Fornberg, B., and Zuev, J.}
\newblock {The Runge phenomenon and spatially variable shape parameters in RBF
  interpolation}.
\newblock {\em Comput. Math. with Appl. 54}, 3 (aug 2007), 379--398.

\bibitem{Freidman1977}
{\sc Freidman, J.~H., Bentley, J.~L., and Finkel, R.~A.}
\newblock {An Algorithm for Finding Best Matches in Logarithmic Expected Time}.
\newblock {\em ACM Trans. Math. Softw. 3}, 3 (1977), 209--226.

\bibitem{FreyGeorge2010}
{\sc Frey, P.~J., and George, P.-L.}
\newblock {\em Quadtree-octree Based Methods}.
\newblock ISTE, 2010, pp.~163--199.

\bibitem{Hardin2005}
{\sc Hardin, D.~P., and Saff, E.~B.}
\newblock {Minimal Riesz Energy Point Configurations for Rectifiable
  d-Dimensional Manifolds}.
\newblock {\em Adv. Math. 193}, 1 (2003), 174--204.

\bibitem{Kuipers2006}
{\sc Kuipers, L., and Niederreiter, H.}
\newblock {\em Uniform distribution of sequences}.
\newblock Dover Publications, Mineola, N.Y, 2006.

\bibitem{Landkof1972}
{\sc Landkof, N.~S.}
\newblock {\em Foundations of modern potential theory}, vol.~180.
\newblock Springer, 1972.

\bibitem{Fornberg2013a}
{\sc Larsson, E., Lehto, E., Heryudono, A., and Fornberg, B.}
\newblock {Stable Computation of Differentiation Matrices and Scattered Node
  Stencils Based on Gaussian Radial Basis Functions}.
\newblock {\em SIAM J. Sci. Comput. 35}, 4 (jan 2013), A2096--A2119.

\bibitem{Lemieux2009}
{\sc Lemieux, C.}
\newblock {\em {Monte Carlo and Quasi-Monte Carlo Sampling}}.
\newblock Springer Series in Statistics. Springer New York, New York, NY, 2009.

\bibitem{Link2012}
{\sc Link, W.~A., and Eaton, M.~J.}
\newblock {On thinning of chains in MCMC}.
\newblock {\em Methods Ecol. Evol. 3}, 1 (feb 2012), 112--115.

\bibitem{Micchelli1986}
{\sc Micchelli, C.~A.}
\newblock {Interpolation of scattered data: Distance matrices and conditionally
  positive definite functions}.
\newblock {\em Constr. Approx. 2}, 1 (1986), 11--22.

\bibitem{Persson2004a}
{\sc Persson, P.-O., and Strang, G.}
\newblock {A Simple Mesh Generator in MATLAB}.
\newblock {\em SIAM Rev. 46}, 2 (jan 2004), 329--345.

\bibitem{Powell2004}
{\sc Powell, M. J.~D.}
\newblock {Five Lectures on Radial Basis Functions}.
\newblock {\em IMM Lect.}, December 2004 (2004), 27.

\bibitem{Preparata1985}
{\sc Preparata, F.~P., and Shamos, M.~I.}
\newblock {\em Computational Geometry}.
\newblock Springer New York, 1985.

\bibitem{Procopiuc_2003}
{\sc Procopiuc, O., Agarwal, P.~K., Arge, L., and Vitter, J.~S.}
\newblock Bkd-tree: A dynamic scalable kd-tree.
\newblock In {\em Advances in Spatial and Temporal Databases}. Springer Berlin
  Heidelberg, 2003, pp.~46--65.

\bibitem{Schoenberg1938}
{\sc Schoenberg, I.~J.}
\newblock {Metric Spaces and Completely Monotone Functions}.
\newblock {\em Ann. Math. 39}, 4 (1938), 811.

\bibitem{Shankar2014}
{\sc Shankar, V., Wright, G.~B., Kirby, R.~M., and Fogelson, A.~L.}
\newblock A radial basis function ({RBF})-finite difference ({FD}) method for
  diffusion and reaction{\textendash}diffusion equations on surfaces.
\newblock {\em Journal of Scientific Computing 63}, 3 (2014), 745--768.

\bibitem{Shimada1995}
{\sc Shimada, K., and Gossard, D.~C.}
\newblock {Bubble mesh}.
\newblock In {\em Proc. third ACM Symp. Solid Model. Appl. - SMA '95\/} (New
  York, New York, USA, 1995), ACM Press, pp.~409--419.

\bibitem{Varma2004}
{\sc Varma, U.~M., Rao, S. V.~R., and Deshpande, S.~M.}
\newblock {Point distribution generation using hierarchical data structures}.
\newblock In {\em Proc. ECCOMAS 2004\/} (Jyv\"askyl\"a, 2004).

\bibitem{VMFF17matlab}
{\sc Vlasiuk, O., and Michaels, T.}
\newblock Boxed lattices and {Riesz} minimizers for {RBF} computations.
\newblock \url{https://github.com/OVlasiuk/3dRBFnodes.git}, 2017.
\newblock [Online; accessed 12-October-2017].

\bibitem{Wendland2004}
{\sc Wendland, H.}
\newblock {\em Scattered Data Approximation}.
\newblock Cambridge University Press, 2004.

\bibitem{Zhou2008}
{\sc Zhou, K., Hou, Q., Wang, R., and Guo, B.}
\newblock {Real-time KD-tree construction on graphics hardware}.
\newblock {\em ACM Trans. Graph. 27}, 5 (dec 2008), 1.

\end{thebibliography}

\end{document}